\documentclass[a4paper,1&pt,twoside]{article}
\usepackage[pdftex,colorlinks,bookmarksnumbered,bookmarksopen,citecolor=red,urlcolor=red]{hyperref}
\usepackage{graphicx}\usepackage{a4wide}
\usepackage{amssymb,amsmath,amsthm}
\usepackage{siunitx}
\usepackage{subfigure,color,colordvi,graphicx,xcolor}
\usepackage[only,llbracket,rrbracket,interleave]{stmaryrd}
\usepackage[sort,compress]{cite}
\usepackage{enumitem}
\usepackage{sidecap, caption}
\usepackage{multirow}

\newcommand{\abs}[1]{\lvert#1\rvert}

\newcommand{\bsym}[1]{\boldsymbol{#1}}
\newcommand{\vms}[1]{\bsym{#1}}

\def\undertilde#1{\mathord{\vtop{\ialign{##\crcr
$\hfil\displaystyle{#1}\hfil$\crcr\noalign{\kern1.5pt\nointerlineskip}
$\hfil\widetilde{}\hfil$\crcr\noalign{\kern1.5pt}}}}}
\newcommand{\xx}{\mathbf{x}}

\newcommand{\ww}{\mathbf{w}}
\newcommand{\wwS}{\ww_{\star}}
\newcommand{\zz}{\mathbf{z}}
\newcommand{\zzS}{\zz_{\star}}

\newcommand{\bS}{{\mathbf{S}}}

\newcommand{\bmu}{\vms{\mu}}
\newcommand{\bmuS}{\vms{\mu}_{\star}}
\newcommand\RR{\leavevmode\hbox{$\rm I\!R$}}
\newcommand{\nS}{\textnormal{\ensuremath{\texttt{n}_{\texttt{s}}}}}
\newcommand{\nP}{\textnormal{\ensuremath{\texttt{n}_{\texttt{p}}}}}
\newcommand{\nO}{{\texttt{n}_{\texttt{o}}}}
\newcommand{\nVar}{\textnormal{\ensuremath{\texttt{n}_{\texttt{var}}}}}
\newcommand{\Pm}{\mathbb{P}_{m}}
\newcommand{\bM}{\mathbf{M}}

\newcommand{\bP}{\mathbf{P}}

\newcommand{\ba}{\mathbf{a}}
\newcommand{\bb}{\mathbf{b}}
\newcommand{\balpha}{\vms{\alpha}}

\usepackage{xparse}
\NewDocumentCommand{\Smu}{o o}{%
  \mathrm{S}_{\texttt{\IfValueTF{#1}{#1}{}}}(\mu_{\IfValueTF{#2}{#2}{}})
}
\NewDocumentCommand{\Sae}{o}{%
  \mathrm{S}_{\texttt{AE}}(\mathbf{x}_{\IfValueTF{#1}{#1}{}})
}
\usepackage{xparse}
\NewDocumentCommand{\bSmu}{o o}{%
  \bS_{\texttt{\IfValueTF{#1}{#1}{}}}(\bmu_{\IfValueTF{#2}{#2}{}})
}

\newcommand{\bm}[1]{\text{\boldmath $#1$\unboldmath}}

\newcommand\ndofx{\nO}
\newcommand\bK{\bm{K}}
\newcommand\bff{\mathbf{f}}
\newcommand\bX{\bm{X}}
\newcommand\bU{\bm{U}}
\newcommand\bUt{\widetilde{\bm{U}}}
\newcommand\bVt{\widetilde{\bm{V}}}
\newcommand\bVtT{\widetilde{\bm{V}}^{\top}}

\newcommand\bZ{\bm{Z}}

\newcommand\bSigma{\bm{\Sigma}}
\newcommand{\nk}{{\texttt{n}_{\texttt{POD}}}}
\newcommand{\bTheta}{\bm{\Theta}}

\newcommand{\eltwo}{\mathcal{L}_2}
\newcommand{\elone}{\mathcal{L}_1}
\newcommand{\XX}{\mathbb{X}}
\newcommand{\Xpgd}{\mathbb{X}_{\texttt{PGD}}}
\newcommand{\Kpgd}{\mathbb{K}_{\texttt{PGD}}}
\newcommand{\Fpgd}{\mathbb{F}_{\texttt{PGD}}}
\newcommand{\bPsi}{\bm{\Psi}}
\newcommand{\fX}{\mathbf{f}_{\mathrm{x}}}
\newcommand{\niter}{\texttt{n}_{\texttt{i}}}
\newcommand{\bs}{{\mathbf{s}}}
\newcommand{\rX}{\mathbf{r}_{\mathrm{x}}}

\newcommand{\nlay}{\texttt{n}_{\texttt{l}}}
\newcommand{\nnJ}[1]{\texttt{n}_{\texttt{n}}^{#1}}
\newcommand{\nlayPrev}{\texttt{n}_{\texttt{n}}^{k{-}1}}

\newcommand{\Lmae}{\mathcal{L}_{\!\texttt{MAE}}}
\newcommand{\Lmse}{\mathcal{L}_{\!\texttt{MSE}}}

\DeclareMathOperator*{\argmin}{arg\,min}

\newenvironment{keywords}{\begin{quote}\emph{\textbf{Keywords:}}}{\end{quote}}
\newtheorem{remark}{Remark}

\RequirePackage{algorithm}[0.1] %Following SIAM cls
\usepackage{algorithmic} %Algorithm style, could alternatively use algpseudocode 

\graphicspath{{figsPDF/}{figsPNG/}{figs/}}

\usepackage{hyphenat}
\usepackage{fancyhdr}

\begin{document}
%
%__________________________________________________________________________________
\fancypagestyle{plain}{%
  \renewcommand{\headrulewidth}{0pt}%
  \fancyhead[R]{}
  \fancyhead[L]{}%
%  \fancyfoot[C]{\footnotesize Page \thepage\ of \pageref{LastPage}}%
}
%__________________________________________________________________________________

%==========================================================================
\title{Surrogates for Physics-based and Data-driven \\ Modelling of Parametric Systems: \\ Review and New Perspectives}

\author{
\renewcommand{\thefootnote}{\arabic{footnote}}
			  Matteo Giacomini\footnotemark[1]\textsuperscript{ \ ,}\footnotemark[2] \  and
			  Pedro D\'{i}ez\footnotemark[1]\textsuperscript{ \ ,}\footnotemark[2]\textsuperscript{ \ ,}*
}

\date{}
%________________________________________________________________________
\maketitle

\renewcommand{\thefootnote}{\arabic{footnote}}

\footnotetext[1]{Laboratori de C\`alcul Numeric (LaC\`aN), E.T.S. de Ingenier\'ia de Caminos, Canales y Puertos, Universitat Polit\`ecnica de Catalunya - BarcelonaTech (UPC), Barcelona, Spain.}
\footnotetext[2]{Centre Internacional de M\`etodes Num\`erics en Enginyeria (CIMNE), Barcelona, Spain. \\
$^*$ Corresponding author: Pedro D\'{i}ez \textit{E-mail:} \texttt{pedro.diez@upc.edu}
}

%________________________________________________________________________
\begin{abstract}
Surrogate models provide compact relations between user-defined input parameters and output quantities of interest, enabling the efficient evaluation of complex parametric systems in many-query settings. Such capabilities are essential in a wide range of applications, including optimisation, control, data assimilation, uncertainty quantification, and emerging digital twin technologies in various fields such as manufacturing, personalised healthcare, smart cities, and sustainability.
This article reviews established methodologies for constructing surrogate models exploiting either knowledge of the governing laws and the dynamical structure of the system (\emph{physics-based}) or experimental observations (\emph{data-driven}), as well as hybrid approaches combining these two paradigms.
By revisiting the design of a surrogate model as a functional approximation problem, existing methodologies are reviewed in terms of the choice of (i) a reduced basis and (ii) a suitable approximation criterion.
The paper reviews methodologies pertaining to the field of Scientific Machine Learning, and it aims at synthesising established knowledge, recent advances, and new perspectives on: dimensionality reduction, physics-based, and data-driven surrogate modelling based on proper orthogonal decomposition, proper generalised decomposition, and artificial neural networks; multi-fidelity methods to exploit information from sources with different fidelities; adaptive sampling, enrichment, and data augmentation techniques to enhance the quality of surrogate models.
\end{abstract}

%________________________________________________________________________

%---------------------------------------------------------------------
\section{Introduction}
\label{sec:Intro}
%---------------------------------------------------------------------

Complex engineering systems are generally characterised by a set of user-defined and/or problem-dependent parameters. In realistic scenarios, these systems feature a potentially large number of parameters, often uncertain.
Stemming from the seminal work~\cite{Box-BW-51} by Box and Wilson, response surfaces have been widely employed in computational science and engineering to provide compact, interpretable, and explainable relations between input parameters and output quantities of interest of the system under analysis.
This notion naturally extends to surrogate models of parametric problems, where a low-dimensional (i.e., compact) representation of high-dimensional data is sought.

This manuscript provides a review of established techniques to construct physics-based, data-driven, and hybrid surrogate models. %revisiting such a task as a functional approximation problem in high-dimensional spaces. 
This task is revisited within the framework of a functional approximation problem in high-dimensional spaces. 
This setting enables the examination of existing solutions in the literature in view of two crucial components: 
\begin{itemize}
\item definition of the approximation criterion for functional representation; 
\item construction of a (reduced) basis for the functional dependence of the quantities of interest on the parameters.
\end{itemize}
This perspective is at the basis of the structure of the present review, with Section~\ref{sec:FuncApprox} focusing on approximation strategies, and Sections~\ref{sec:POD},~\ref{sec:PGD}, and~\ref{sec:NN} describing different approaches to construct the reduced basis.
The practical deployment of such techniques to devise low-dimensional parametric representations in realistic, many-query, and data-scarce settings depends on the ability to handle heterogeneous data (i.e., models of different fidelities) efficiently, that is minimising the amount of required information.
For this reason, multi-fidelity modelling (Section~\ref{sec:MultiFidelity}) and adaptive dataset sampling (Section~\ref{sec:Sampling}) are discussed in dedicated sections of the manuscript.
These topics are not specific to a particular class of techniques, but act as cross-cutting enablers enhancing accuracy, robustness, and efficiency across physics-based, data-driven, and hybrid approaches to surrogate modelling.
A detailed presentation of the structure of the manuscript is reported in Section~\ref{sec:Structure}.

Several comprehensive reviews on reduced-order and surrogate modelling for parametric systems are available in the literature, from general monographs~\cite{Benner-BGQRSS-21-v1,Benner-BGQRSS-21-v2,Benner-BGQRSS-21-v3} to reviews of specific topics such as dimensionality reduction~\cite{Hou-HB-22}, projection-based methods for dynamical systems~\cite{Benner-BGW-15}, and reduced-order models as input-output control systems~\cite{Benner-BBF-14}.
The present work complements these contributions by adopting a unifying view on the construction of surrogate models within a functional approximation perspective.
This allows to present an overview of \emph{general purpose} techniques for surrogate modelling of parametric systems, spanning physics-based, data-driven, and hybrid modelling approaches.
In particular, it focuses on the methodological aspects for the development of efficient surrogate models, seamlessly analysing solutions from different fields.
Indeed, many relevant contributions have been made on such topics by different academic communities, from applied mathematics and computational science to domain-specific engineering disciplines.
In this context, the very same notation used for classification is not always unique and consistent across different communities.
In this work, we denote by \emph{ROM} a model solving the governing equations of the system under analysis in a reduced space (e.g., via projection methods), whereas we employ the notion of \emph{surrogate} for a more general class of models providing a functional approximation of a quantity of interest.
This classification does not restrict surrogates to purely data-driven approaches, allowing also for hybrid techniques blending data and physical knowledge.
Finally, it is worth noting that a thorough and fair comparison of existing methodologies is extremely challenging (if even possible) since each method is expected to present advantages and disadvantages depending on the specific applications selected for benchmarking.
Indeed, many surrogate solutions tailored to given physical and engineering problems have been proposed in the literature and readers interested in domain-specific surveys for computational fluid and structural mechanics are referred to~\cite{Rozza-RSB-22,Samadian-SMD-25}.

Following the general purpose rationale mentioned above, the surrogate modelling techniques reviewed in this work are assessed from the methodological viewpoint rather than through application-specific benchmarks.
Throughout the manuscript, the discussion is guided by a set of evaluation criteria reflecting key methodological aspects of surrogate construction. These include the degree of intrusiveness of the low-dimensional representation with respect to the high-fidelity model, the nature of the approximation space (e.g., linear versus nonlinear), the amount of training data required, the robustness with respect to increasing the dimensionality of the system and the parametric variability, and the interpretability and explainability of the resulting surrogate models.
In addition, practical considerations on the performance of the reviewed approaches are considered, including: the computational cost of the offline training stage, the efficiency of the online evaluations, the scalability to many-query settings, and the suitability of the resulting surrogates for real-time analysis, optimisation, or decision-support applications involving multiple quantities of interest.
These criteria are not intended to establish a quantitative ranking of the presented methods. Rather, they provide a common analytical perspective to position the different methodologies within a unified framework, allowing their respective strengths, limitations, and methodological trade-offs to be discussed in a coherent and structured manner. 

In the remainder of this section,  the functional approximation problem is presented and the main concepts involved in the construction of a (parametric) surrogate model are recalled in order to lay the foundations of the following review.

%---------------------------------------------------------------------
\subsection{Parametric full-order models}
\label{sec:FOM}
%---------------------------------------------------------------------

It is common to denote as \emph{full-order} the model describing the behaviour of the system. 
\emph{Physics-based} full-order models rely on a mathematical formulation of the problem by systems of ordinary and/or partial differential equations. In this context, each instance of a full-order model entails the discretisation of such equations, resulting in (possibly large) systems of linear or nonlinear algebraic equations.
\emph{Data-driven} full-order models can be obtained either from experimental observations or from datasets generated by numerical simulators. In many practical settings, these simulators are treated as \emph{black-box} models where only input-output data are accessible (e.g., commercial software with proprietary source codes). In addition, such datasets can also be obtained from empirical or semi-empirical models whose internal formulation is not derived from first-principles physical laws of the underlying system. In all these cases, the surrogate model is inferred directly from the available data without requiring explicit access to, or exploitation of, the governing equations.

In common engineering settings, full-order models need to be evaluated for different combinations of the input parameters. This typically leads to many queries of the model, populating high-dimensional manifolds with solutions. Relevant applications include optimisation, control, inverse analysis, and uncertainty quantification, to name a few.
To make these problems treatable,  a low-dimensional representation of such high-dimensional solutions (or of quantities of interest depending on them) is to be constructed, thus allowing for efficient, real-time evaluations of the parametric model.

%---------------------------------------------------------------------
\subsection{Surrogates as functional approximations}
\label{sec:Surrogates}
%---------------------------------------------------------------------

Let a system be characterised as a function transforming the input parameters $\bmu=[\mu_{1} \, \mu_{2} \dots \mu_{\nP}]^{\top}$ ranging in the set of admissible values $\mathcal{P} \subset \RR^{\nP}$ into an output vector $\xx=[x_{1} \,x_{2} \dots x_{\nO}]^{\top}$,  where $\nP$ and $\nO$ are the number of parameters and the size of the system output, respectively. 
Thus, the output $\xx$ is seen as a function of $\bmu$ to be approximated by a surrogate $\bS:\RR^{\nP} \rightarrow \RR^{\nO}$ such that 
\begin{equation}\label{eq:FuncApp}
\xx(\bmu)\approx \bS(\bmu). 
\end{equation}

Note that equation~\eqref{eq:FuncApp} is the general form of a functional approximation problem, where $\bS$ is the sought low-dimensional representation of the map transforming the input parameters $\bmu$ into the high-dimensional data $\xx(\bmu)$.
In this context, $\xx$ could represent the solution of some partial differential equation obtained from a physics-based full-order model, as well as a quantity of interest, e.g. a force measured experimentally without any knowledge of the underlying mathematical model.

As is typical in functional approximation problems, four aspects are to be taken into account.
\begin{itemize}
\item \textbf{Knowledge of the true system -} Available information (if any) on the function $\xx(\bmu)$ linking inputs and outputs, that is,  any a priori knowledge on the behaviour of the true system, stemming from either experimental observations or mathematical formulations.
\item \textbf{Available data -} Information about observable data (i.e. inputs and outputs) of the function $\xx(\bmu)$. Data often consist of a set of discrete points, that is pairs $(\bmu_{i},\xx_{i})$ for $i=1,\ldots, \nS$,  $\nS$ being the number of samples such that $\xx_{i}:=\xx(\bmu_{i})$.
\item \textbf{Approximation criterion -} Choice of the strategy to establish \emph{how} to fit available data. The approximation criterion specifies the rule employed to establish the \emph{best fit}, that is  what ``$\approx$'' means in equation~\eqref{eq:FuncApp}. In practice, this criterion provides a suitable numerical measure of the mismatch between the response function $\bS(\bmu_i)$ obtained from the surrogate model and the true data $\xx_i$, for $i=1,\ldots,\nS$. A review of different approximation criteria will be presented in Section~\ref{sec:FuncApprox}.
\item \textbf{Type of approximation -} Selected basis, that is, nature of the function employed to represent $\bS(\bmu)$.  There is a vast choice of function types, such as polynomials, splines,  other functions with analytical expression (e.g. radial basis functions), neural networks or alternative functions to be determined in a constructive manner. In the context of high-dimensional data, it is crucial to select a basis allowing to reduce the dimensionality of the functional representation, identifying a suitable low-rank approximation or reduced-order model (ROM)~\cite{AH-CHRW:17}. Different strategies to construct such a reduced basis are detailed in Sections~\ref{sec:POD}, \ref{sec:PGD}, and~\ref{sec:NN}. 
\end{itemize}

Note that the four points mentioned above are strongly correlated, and available data influences the choice of the functional type and the approximation criterion.

%---------------------------------------------------------------------
\subsection{Verification, validation, and uncertainty quantification}
\label{sec:VVUQ}
%---------------------------------------------------------------------

Surrogate models are often employed to accelerate \emph{outer-loop} applications such as optimisation, control, and inverse analysis, to name a few.
This entails repeated evaluations of quantities of interest under parametric variability, hence assessing the trade-off between efficiency and credibility of such computations is of paramount importance.
Evaluating and controlling the errors due to the mismatch of the model with respect to reality and to the employed numerical discretisations (i.e., verification and validation~\cite{Babuska-Oden-04}) play a crucial role to guarantee reliability in decision protocols involving realistic problems.
These challenges entail the need for accurate, efficient, and reliable error estimates to be coupled with model updating strategies and adaptivity procedures to improve the overall quality of the surrogate models.

In addition, accounting for data uncertainties~\cite{UQ-book} is crucial in the construction of \emph{credible} surrogate models.
On the one hand, user-defined parameters are frequently affected by noise and are known up to certain measuring precisions. On the other hand, it is common for operating conditions of complex systems to be uncertain.
Hence, quantifying the uncertainty on system outputs is crucial to perform inference and predict the range of reliability of the quantities of interest in the context of high-dimensional, many-queries, complex systems.

In this context, contributions related to verification and validation (V\&V) and uncertainty quantification (UQ) naturally constitute key aspects in the literature on surrogate models for parametric systems. 
Several critical performance aspects of surrogate models naturally arise in this setting, including extrapolation capabilities, error estimation, confidence bounds, robustness with respect to parametric variations, and stability or divergence properties.
The assessment of such properties is highly problem-dependent and typically requires a careful interplay between the surrogate construction methodology, the characteristics of the underlying physical system, and the adopted V\&V and UQ frameworks.
Whilst these aspects are extremely relevant in the context of reduced and surrogate models, this work focuses on the methodological foundations for the construction of general purpose surrogate models for parametric systems, encompassing physics-based, data-driven, and hybrid approaches.
A comprehensive review of V\&V, UQ, and adaptive error control strategies for surrogate models is beyond the scope of the present work and interested readers are referred, e.g., to the review article~\cite{Gunzburger-PWG-18}.

%---------------------------------------------------------------------
\subsection{Structure of the manuscript}
\label{sec:Structure}
%---------------------------------------------------------------------

The manuscript is organised as follows.
Section~\ref{sec:FuncApprox} revisits the formulation of different approximation criteria for functional approximation in dimension $1$, from polynomial interpolation to least-squares approximation, as well as more advanced techniques such as radial basis functions, Kriging, moving least-squares, and universal approximation via artificial neural networks.
The following three sections address the construction of the reduced basis by means of proper orthogonal decomposition (Section~\ref{sec:POD}), proper generalised decomposition (Section~\ref{sec:PGD}), and neural networks (Section~\ref{sec:NN}). For each methodology, the corresponding section reviews existing contributions in the literature, covering dimensionality reduction and construction of purely data-driven and physics-based surrogates. 
Section~\ref{sec:MultiFidelity} is devoted to multi-fidelity models, reviewing existing methodologies to fuse data of different fidelities and sources.
The challenges of sampling, database enrichment, and data augmentation are discussed in Section~\ref{sec:Sampling}.
Finally, Section~\ref{sec:Conclusion} summarises the discussion presented in the article and presents open problems, current trends, and future perspectives of the research on surrogate models constituting key building blocks to enable efficient, scalable, and trustworthy digital twin technologies.

%---------------------------------------------------------------------
\section{Revisiting functional approximation criteria}
\label{sec:FuncApprox}
%---------------------------------------------------------------------

Although approximation criteria for surrogate modelling are ultimately applied within a chosen reduced representation, they can be formulated independently of the specific construction of the reduced space.
For this reason, this section first recalls the fundamentals of functional approximation, both based on interpolation and regression. Sections~\ref{sec:POD},~\ref{sec:PGD}, and~\ref{sec:NN} then introduce three different methodologies to devise suitable reduced bases which, combined with the approximation strategies discussed here, enable the construction of parametric surrogate models.

Let $\nP=1$ and $\nO=1$. The corresponding approximation problem involves determining a surrogate $\Smu$ of the scalar function $x(\mu)$, starting from a set of $\nS$ data pairs $(\mu_i,x_i)$,  with $x_i = x(\mu_i)$,  for $i=1,\ldots, \nS$.

%---------------------------------------------------------------------
\subsection{Polynomial interpolation}
\label{sec:Poly}
%---------------------------------------------------------------------

Polynomial interpolation consists of determining as response function a polynomial $\Smu[L]$ of degree $\nS-1$ that matches the data pairs $(\mu_i,x_i)$, that is, $\Smu[L][i]=x_i$,  for $i=1,\ldots, \nS$. 

This requires setting the degree of the polynomial equal to the number of samples $\nS$ minus one.
If the sampling points are different (i.e. , $\mu_i \neq \mu_j$ for $i\neq j$) then there is a unique polynomial function $\Smu[L]$ of degree $\nS-1$ fitting the data. 

The polynomial interpolation response is thus computed as linear combination of a set of polynomial basis functions, that is, 
\begin{equation}\label{eq:LagrangeInterp}
\Smu[L] := \sum_{i=1}^{\nS}  x_i \, L_i(\mu) ,
\end{equation}
where $L_i(\mu)$ denotes the $i$-th Lagrange polynomial function defined as
\begin{equation}\label{eq:LagrangePoly}
L_i(\mu) := \frac{\displaystyle\prod_{j \neq i} \mu-\mu_j}{\displaystyle\prod_{j \neq i} \mu_i-\mu_j} .
\end{equation}

It is straightforward to extend this idea to multiple dimensions, that is,  $\nP > 1$. 
The definition of the corresponding polynomial basis generalising \eqref{eq:LagrangePoly} has to be devised to obtain $L_i$ such that $L_i(\bmu_j) = \delta_{ij}$, with $\delta_{ij}$ being the Kronecker delta. 
This is straightforward if samples are available in a Cartesian lattice, that is, if the parametric space is discretised using a grid obtained as the tensorial product of $\nP$ one-dimensional discretisations, one for each parametric dimension. 

Note that the computational effort and data complexity increases exponentially with $\nP$, awakening the so-called \emph{curse of dimensionality}.
Strategies to alleviate this computational burden, especially in high-dimensional settings, have been proposed by suitably selecting a reduced number of sampling points in the Cartesian grid, e.g., relying on the sparse grids framework~\cite{SparseGrids-04}.
Exploiting these points,  the sparse grid stochastic collocation method~\cite{Nobile-SIAMJNA-08} allows to significantly reduce the computational cost, while maintaining the ability to fit a polynomial response function.

\begin{remark}[Spectral approximations]
Alternative choices of basis functions are also possible, e.g., replacing~\eqref{eq:LagrangeInterp} by means of orthogonal polynomials to construct spectral approximations~\cite{Solin-SSD-book-03}.
When noisy evaluations are to be accounted for, this interpolation procedure generalises to a regression problem.
In this context, polynomial chaos expansion (PCE)~\cite{Ghanem-GS-03,Karniadakis-XK-02,Karniadakis-XK-03} can be viewed as a form of polynomial interpolation using an orthogonal basis over a random input space.
In particular, PCE coefficients are estimated in a least-squares sense, projecting the mean response onto the orthogonal polynomial basis while filtering out noise.
This framework has been extensively employed to construct surrogate models, especially for uncertainty quantification, see the foundational works~\cite{Xiu-10,LeMaitre-MK-10,Sudret-BS-10,Sudret-BS-11} and the recent review~\cite{Sudret-LMS-21}.
\end{remark}

%---------------------------------------------------------------------
\subsection{Non-polynomial interpolation}
\label{sec:NonPoly}
%---------------------------------------------------------------------

Interpolation-based functional approximations are also devised letting users select alternative types of basis functions.
For example, trigonometric functions could be used instead of polynomials, assuming that the proposed functional basis is also complete. 
Indeed,  as Taylor expansion provides an approximation of any given function in terms of a sum of polynomials, Fourier expansion employs a sum of trigonometric terms for the same task.

Alternative techniques in which the type of functional approximation (i.e. the functional space $\Smu$ belongs to) is much more flexible have also been proposed in the literature. 
In this section, radial basis functions (RBFs) and Gaussian process regression (GPR) are reviewed.

The paradigm of RBF interpolation consists of building an approximation $\Smu[RBF]$ of the form
\begin{equation}\label{eq:RBFinterp}
\Smu[RBF] := \sum_{i=1}^{\nS} w_i  \phi(\vert \mu - \mu_i \vert)
\end{equation}
where $w_i$ denotes the $i$-th weight and $\phi$ is a given radial basis function that depends on the distance between the point $\mu$, where the function has to be evaluated, and the sampling points $\mu_i$, for $i=1,\ldots, \nS$. 
Typically, the radial basis function $\phi(r)$ is defined such that $\phi(0)=1$ and $\phi(r)$ decreases as $r$ increases, for instance, the expression of the Gaussian RBF is $\phi(r)=e^{-\varepsilon^2 r^2}$~\cite{Powell-RBF-87}.

Given the set of $\nS$ sample points $(\mu_i,x_i)$, the weights $w_i$ are computed to enforce the interpolation property, that is, $\Smu[RBF][i]=x_i$,  for $i=1,\ldots, \nS$. 
This results in solving the linear system of equations 
\begin{equation}\label{eq:RBFsystem}
\bm{A}
\begin{bmatrix}
w_{1}\\
w_{2}\\
\vdots \\
w_{\nS}
\end{bmatrix}
=
\begin{bmatrix}
x_{1}\\
x_{2}\\
\vdots \\
x_{\nS}
\end{bmatrix} ,
\end{equation}
with $\bm{A}$ being the matrix of dimension $\nS \times \nS$ with components $[\bm{A}]_{ij}$ defined as
\begin{equation}\label{eq:RBFImatrix}
[\bm{A}]_{ij} := \phi( \vert \mu_{i} {-} \mu_{j} \vert) .
\end{equation}
Note that the diagonal terms of $\bm{A}$ are all equal to 1 and that the non-diagonal terms are lower than one and often negligible if the distance between the associated sampling points is large enough.

Under the assumption of noise-free data, Gaussian process regression (GPR), also known as Kriging, is an interpolation technique in which the resulting functional approximation is defined implicitly through kernel evaluations and training data, rather than through an explicit parametric analytical expression.
It is worth noting that the term \emph{regression} is to be intended here in a broad statistical sense, referring to the task of learning the input-output relation to make predictions for unseen data.
The Kriging rationale assumes that the value of the response function $\Smu[GPR]$ at every point is a random variable with a Gaussian (or normal) distribution. 
It is also assumed that the joint distribution of the variables corresponding to two different points is a bivariate normal and that the associated covariance depends only on the distance $r$ between the points. This covariance is thus given by a kernel function $\kappa(r)$ that plays a similar role as the radial basis function $\phi$ described above.

Following the idea of Bayesian updating, the approximation $\Smu[GPR]$ is found by introducing the information brought by the samples as an update of a prior estimate, that may be as uninformative as the user wishes.
The matrix of variances and covariances of the samples, denoted by $\bm{C}$,  has an expression very similar to \eqref{eq:RBFImatrix}. In particular, for  $i=1,\ldots,\nS$ and for $j=1,\ldots,\nS$, the component $[\bm{C}]_{ij}$ of the matrix is defined as
\begin{equation}\label{eq:KrigingMatrix}
[\bm{C}]_{ij} := \kappa( \vert \mu_{i} {-} \mu_{j} \vert) .
\end{equation}
Similarly, at each point $\mu_\star$ where the Kriging interpolation is sought,  the vector $\bm{\tau}$ defining the spatial correlation with the $\nS$ samples is computed as
\begin{equation}\label{eq:KrigingNewEst}
\bm{\tau} =
\begin{bmatrix}
\kappa( \vert \mu_{1}- \mu_\star \vert)\\
\kappa( \vert \mu_{2}- \mu_\star \vert)\\
\vdots \\
\kappa( \vert \mu_{\nS}- \mu_\star \vert)
\end{bmatrix} .
\end{equation}

\begin{remark}[Variogram]
Instead of directly providing the covariance using the kernel function $\kappa$, a common alternative relies on computing the variogram~\cite{Oliver-Webster-15}. 
The variogram quantifies the dissimilarity (or variability) between data points to model underlying spatial relations.
In practice, the variogram is constructed by calculating, for the pairs of points separated by a distance $r$, and for different values of $r$,
\begin{subequations}
\begin{equation}\label{eq:variogram}
\gamma(r) = \frac{1}{2\nVar(r)} \sum_{i=1}^{\nS} \sum_{j=1}^{\nS} \beta_{ij}(r) (x_i-x_j)^2 ,
\end{equation}
where the summation in~\eqref{eq:variogram} is restricted only to the $\nVar(r)$ sets of pairs $(i,j)$ whose positions are approximately by a distance $r$, as identified by the coefficients $\beta_{ij}(r)$ such that
\begin{equation}\label{eq:variogramCoef}
\beta_{ij}(r) = 
\begin{cases}
1 , & \text{if $| \mu_i {-} \mu_j | \in [r-\eta/2,r+\eta/2]$}, \\
0 , & \text{otherwise}
\end{cases}
\end{equation}
with $\eta {>} 0$ being a user-defined tolerance to approximate the distance $r$.
The function in~\eqref{eq:variogram} thus describes how variability of data increases with spatial separation. A theoretical variogram model (e.g., spherical, exponential, or Gaussian) is then fitted to the empirical values to obtain a continuous representation of spatial variability. 
Under the assumption of second-order stationarity, the variogram can be converted into the corresponding covariance function
\begin{equation}\label{eq:variogramCovariance}
\kappa(r) = \kappa(0) - \gamma(r) ,
\end{equation}
\end{subequations}
which is subsequently used in the Kriging procedure to perform spatial interpolation.
\end{remark}

Hence, in the absence of noise, the Kriging interpolated value $\Smu[GPR][\star]$ is obtained as a weighted sum of data $x_i$, with weights coming from the solution of the linear system of equations $\bm{C} \bm{\lambda} = \bm{\tau}$, that is,
\begin{equation}\label{eq:interpGPR}
\Smu[GPR][\star] := \sum_{i=1}^{\nS}  x_i \lambda_i  ,
\end{equation}
where the weights $\lambda_i$ are calculated to minimise the variance of the prediction while ensuring that the interpolated value is unbiased and the symmetry of matrix $\bm{C}$ is exploited.

Whilst in the case of noise-free data Gaussian process regression provides an interpolatory response surface, in the presence of noise GPR incorporates a noise variance term, producing a smooth surrogate which no longer interpolates the observed data.
The regularity of the resulting surrogate depends on the choice of the covariance kernel, with different kernels inducing different degrees of smoothness of the predicted response.
Moreover, when the noise level is not constant over the parameter space, techniques such as heteroscedastic Gaussian process regression~\cite{Canu-LSC-05} or stochastic Kriging~\cite{Staum-ANS-10} should be considered.
Taking the prior estimate of $\Smu[GPR][\star]$ as a normal distribution with mean $\theta_{\text{prior}}$ and variance  $\sigma^{2}_{\text{prior}}$, the updated (or posterior) mean and variance estimates accounting for the information at hand read
\begin{subequations}
 \begin{align}
\theta_{\text{posterior}} & = \theta_{\text{prior}} +\bm{\tau}^{\top} \bm{C}^{-1}(\xx-\bm{\theta}_{\text{prior}} ) , \\
\sigma^{2}_{\text{posterior}} & = \sigma^{2}_{\text{prior}} - \bm{\tau}^{\top} \bm{C}^{-1}\bm{\tau} ,
\end{align}
\end{subequations}
where $\xx$ is the vector of samples $\xx=[x_{1} \, x_{2} \dots x_{\nS}]^{\top}$ and $\bm{\theta}_{\text{prior}} $ is the vector of the prior estimated mean in the sampling points.
Hence,  the posterior mean provides the best estimate at a new sample $\mu_\star$ given the observed data $(\mu_i,x_i)$ for $i=1,\ldots,\nS$ and the spatial correlation $\bm{\tau}$ of the model.  

Note that, from an algorithmic perspective, both RBF and Kriging interpolation can be extended to high-dimensional parameter spaces ($\nP \gg 1$), since their computational complexity depends on the number of samples $\nS$ rather than on the number of parameters $\nP$. 
Of course, the number of samples required to construct an accurate surrogate model generally increases with the number of parameters $\nP$, whence Section~\ref{sec:Sampling} provides a detailed discussion of efficient sampling strategies.

%---------------------------------------------------------------------
\subsection{Least-squares approximation}
\label{sec:LS}
%---------------------------------------------------------------------

Let $\Pm$ be the space of polynomials up to degree $m$.
The least-square approximation $\Smu[LS]$ is sought in the linear approximation space $\Pm$ and has the form
\begin{equation}\label{eq:LSapprox1D}
\Smu[LS] :=\sum_{k=0}^m a_k P_k(\mu) ,
\end{equation}
where $\{P_0 \, P_1 \dots P_m\}$ denotes a basis of the space $\Pm$.  Typically these bases are such that $P_k$ is  a polynomial of degree $k$; for instance, in the case of a monomial basis, it holds $P_k(\mu) = \mu^k$.

Let $\langle \cdot,\cdot \rangle$ denote a scalar product in the parametric space $\mathcal{P}$.
The scalar product can be defined as the continuous form
\begin{subequations}
\begin{equation}\label{eq:continuousSC}
\langle f,g \rangle=\int_{\mathcal{P}} \alpha(\mu) f(\mu) g(\mu)  \, d\mu ,
\end{equation}
where $\alpha(\mu)$ represents a suitably defined weighting function, or at the discrete level as
\begin{equation}\label{eq:discreteSC}
\langle f,g \rangle= \sum_{i=1}^{\nS} \alpha_i f(\mu_i) g(\mu_i)  .
\end{equation}
\end{subequations}

The least-squares criterion consists of finding $\Smu[LS]$ in order to minimise the distance
\begin{equation}\label{eq:distLS}
\begin{aligned}
\Vert \Smu[LS] &- x(\mu) \Vert^2 \\
&=\langle \Smu[LS]-x(\mu),\Smu[LS]-x(\mu) \rangle ,
\end{aligned}
\end{equation}
that corresponds to solving the optimisation problem
\begin{equation}\label{eq:minLS}
\Smu[LS] = \arg \min_{\mathrm{S}_\star\in \Pm} \Vert \mathrm{S}_\star(\mu)-x(\mu) \Vert^2 .
\end{equation}
Note that the function $\Smu[LS]$ provides the \emph{best} global approximation for all data in the least-squares sense, but does not guarantee interpolation properties, that is $\Smu[LS][i] \neq x_i$, for any $i=1,\ldots,\nS$.

The least-squares solution is thus obtained by solving the normal equations, that is, the linear system of algebraic equations
\begin{equation}\label{eq:normalEqs}
\bM \ba = \bb ,
\end{equation}
where $\ba=[a_0 \, a_1 \dots a_m]^\top$ is the vector of the unknown coefficients in~\eqref{eq:LSapprox1D},  $\bM$ is the $(m{+}1){\times}(m{+}1)$ matrix such that entry $[\bM]_{k\ell}:=\langle P_k,P_\ell \rangle$,  and $\bb$ is the vector with entries given by $[\bb]_k:=\langle P_k,x \rangle$.
If the basis $\{P_0 \, P_1 \dots P_m\}$ is orthogonal (i.e., $\langle P_k,P_\ell \rangle=0$ for $k\neq \ell$), then the coefficients are such that $a_k=\langle x,P_k \rangle/\langle P_k,P_k \rangle$.
Using the discrete scalar product~\eqref{eq:discreteSC},  the resulting least-squares approximation is given by
\begin{equation}\label{eq:LS}
\Smu[LS] :=\sum_{k=0}^m \left[\frac{\displaystyle\sum_{i=1}^{\nS} \alpha_i x_i P_k(\mu_i)}{\displaystyle\sum_{i=1}^{\nS} \alpha_i P_k(\mu_i) P_k(\mu_i)} \right] \! P_k(\mu) .
\end{equation}

\begin{remark}[Choice of the scalar product]
Each scalar product is associated with a different basis of orthogonal polynomials.
Hence, the choice of the scalar product $\langle \cdot,\cdot \rangle$, in particular the weighting function $\alpha$, sets the approximation properties.
\end{remark}

Alternatively to the least-squares approach, the minimax criterion, which also features an approximation of the form~\eqref{eq:LSapprox1D}, proposes to determine the function by minimising the maximum error in $\mathcal{P}$, leading to the problem
\begin{equation}\label{eq:minimax}
\Smu[minimax] = \arg \min_{\mathrm{S}_\star\in \Pm} \max_{\mu_\star \in \mathcal{P}} \vert \mathrm{S}_\star(\mu_\star) - x(\mu_\star)  \vert .
\end{equation}
Whilst the minimax criterion is attractive as it allows to control the approximation in the worst-case scenario (and not its average as in least-squares),  it results in a nonlinear system of equations to determine the unknown weights in~\eqref{eq:LSapprox1D}.
Nonetheless, selecting $\alpha(\mu)=1/\sqrt{1-\mu^2}$ in~\eqref{eq:continuousSC}, the resulting family of orthogonal polynomials is known as the Chebyshev polynomials, and the least-squares criterion (that requires solving only the linear problem~\eqref{eq:normalEqs}) provides a very good approximation of the minimax solution $\Smu[minimax]$.
Indeed, if $x(\mu)$ is a polynomial of degree $m+1$, the Chebyshev approximation is such that $x(\mu)-\Smu[LS] \propto P_{m+1}(\mu)$. Moreover,  $P_{m+1}(\mu)$ being a function that oscillates $m$ times between $-1$ and $+1$, it follows that $\Smu[LS]$ precisely coincides with the minimax approximation.

\begin{remark}[Optimal sampling points]\label{rmrk:Cheb}
If the information about the function $x(\mu)$ is discrete (i.e., values at different points $\mu_i$, for $i=1,\ldots,\nS$), the discrete form \eqref{eq:discreteSC} of the scalar product has to be adopted to compute $\bb$.
This naturally raises the question of how to optimally select the sampling points.
A classical choice consists in selecting the $\nS$ points as the zeros of the Chebyshev polynomial of degree $\nS{-}1$. In this case, the discrete orthogonality induced by the associated quadrature rule yields polynomial bases that are closely related to the continuous Chebyshev polynomials, and recover their approximation properties up to the degree for which the quadrature is exact.
As a consequence, the resulting approximation exhibits near-minimax behaviour and effectively mitigates the Runge phenomenon, making Chebyshev nodes a \emph{near-optimal} choice for polynomial interpolation from discrete data.
\end{remark}

%---------------------------------------------------------------------
\subsection{Moving least-squares approximation }
\label{sec:MLS}
%---------------------------------------------------------------------

The moving least-squares (MLS) approximation is based on the standard least-squares approach presented above, following the rationale of the smoothed particle hydrodynamics (SPH) method.
For a complete introduction to the subject, interested readers are referred to~\cite{AH-HBFRZA:17}.

The SPH idea from \cite{Lucy-77} consists of approximating a function $x(\mu)$ as
\begin{equation}\label{eq:approxSPH}
\Smu[SPH] := \sum_{i=1}^{\nS}  \underbrace{C_\rho \omega_i \Phi\left( \frac{\mu-\mu_i}{\rho} \right)}_{=:N_i(\mu)} x_i ,
\end{equation}
where $\Phi$ is a window function (also seen as a kernel), $\rho$ is a dilation parameter accounting for the radius of influence of each sampling point $\mu_i$, $\omega_i$ measures the area (or hyper-volume) of influence of the sampling point, and $C_\rho$ is a normalisation coefficient. The resulting basis function $N_i(\mu)$ is a simple scaling of the window function $\Phi$.

MLS follows a similar formulation, exploiting the same basis $\{P_0 \, P_1 \dots P_m\}$ of the space $\Pm$ introduced in Section~\ref{sec:LS} for the standard least-squares approach. 
The idea is taking the basis function 
\begin{equation}\label{eq:basisMLS}
N_i^{\texttt{MLS}} (\mu) := \balpha(\mu)^\top \bP(\mu_i) \Phi\left( \frac{\mu-\mu_i}{\rho} \right) ,
\end{equation}
where $\bP=[ P_0 \, P_1 \dots P_m ]^\top$ gathers the polynomial basis functions up to degree $m$, and $\balpha$ is the vector of unknown functions of  $\mu$ to be determined.
This leads to solving a linear system of equations to compute $\balpha$ \emph{for all values of} $\mu$. That is, in all the points $\mu_i, \, i=1,\ldots,\nS$ where the surrogate is to be computed, it holds
\begin{equation}\label{eq:EqMLS}
\bM(\mu) \balpha(\mu) = \bP(\mu) ,
\end{equation}
with $\bM(\mu)$ being a weighted mass matrix (or Gramm matrix) defined as 
\begin{equation}\label{eq:matMLS}
\bM(\mu):= \sum_{i=1}^{\nS}  \bP(\mu_i)  \bP(\mu_i)^\top  \Phi\left( \frac{\mu-\mu_i}{\rho} \right) .
\end{equation}
Note that the sum in~\eqref{eq:matMLS} can be further restricted to consider exclusively the points \emph{close} to the value $\mu_i$ where the window function is centred. 
In this context, the notion of \emph{vicinity} is to be understood in relation with the dilation parameter $\rho$.

The resulting MLS approximation is thus given by
\begin{equation}\label{eq:MLS}
\Smu[MLS] :=\sum_{i=1}^{\nS} x_i \, N_i^{\texttt{MLS}}(\mu) ,
\end{equation}
and it fulfils the following properties:
\begin{itemize}
\item \textbf{Consistency or reproducibility of order $m$ -} The approximation is exact for $x(\mu)$ being a polynomial of degree up to $m$.
\item \textbf{Adaptivity -} The accuracy of the discretisation is easily enhanced locally by adding new sampling points and adapting the dilation parameter $\rho$ such that it is proportional to the characteristic distance between points (typically denoted as $h$).
\item \textbf{Regularity -}  MLS shape functions $N_i^{\texttt{MLS}} (\mu)$ and the corresponding approximation $\Smu[MLS]$ have the same regularity as the window function $\Phi$.
\item \textbf{Non-interpolative -} Similarly to the standard least-squares regression, the approximation $\Smu[MLS]$  is not interpolative, that is in general $\Smu[MLS][i] \neq x_i$, for  $i=1,\ldots,\nS$.
\end{itemize}

%---------------------------------------------------------------------
\subsection{Universal approximation via artificial neural networks}
\label{sec:introNN}
%---------------------------------------------------------------------

Artificial neural networks or, simply, Neural Networks (NN) are response models linking input and output parameters via a set of user-defined layers of \emph{neurons}. 
Feed-forward NNs naturally provide a flexible framework for functional approximation by means of a nonlinear regression strategy in which the form of the function is implicitly learned from data instead of being defined a priori by the user.
Following the seminal result of the universal approximation theorem~\cite{UniversalThNN-89}, feed-forward nets with at least one hidden layer and suitable nonlinear activation functions can be used to approximate any continuous function to an arbitrary level of accuracy.

\begin{figure*} [h]
	\centering	
	\includegraphics[width=\textwidth]{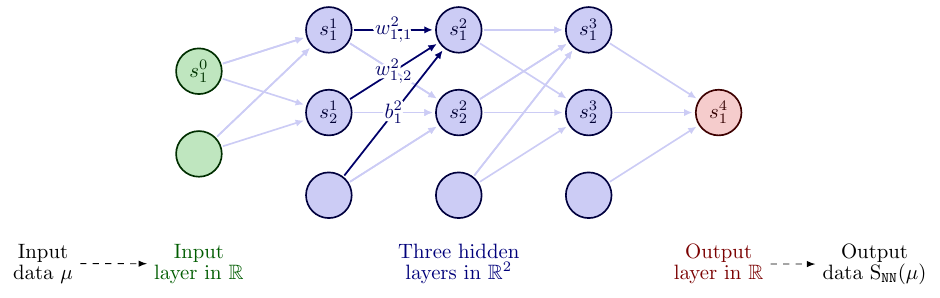}
	\caption{Illustration of a neural network with one input neuron, three hidden layers with two neurons and one output neuron. In each layer, the blank dot is associated with bias.}
	\label{fig:NN}
\end{figure*}

Let $\nlay$ be the number of layers of the net and denote by $\nnJ{k}$ the number of neurons in the $k$-th layer. The state of the $j$-th neuron in the $k$-th layer of the net is represented by $s_j^k$.
In a fully-connected architecture, each node of the $k$-th layer is connected to all the neurons in the $(k{-}1)$-th layer and a weight $w_{j,\ell}^k$ is associated with each connection. Moreover, each neuron carries an intrinsic bias $b_j^k$.
A simple architecture featuring one input parameter in the first layer, three hidden layers with two neurons, and an output layer with one neuron is depicted in Figure~\ref{fig:NN}.

The state of each neuron is defined as
\begin{equation}\label{eq:NNstate}
s_j^k = g\!\left( \sum_{\ell=1}^{\nlayPrev} w_{j,\ell}^k s_\ell^{k-1} + b_j^k  \right) ,
\end{equation}
where $g$ is the so-called \emph{activation function}, a (generally nonlinear) function that transforms the inputs and bias of each neuron into its output state.
Commonly employed activation functions are the rectified linear unit (ReLU), the hyperbolic tangent, and the sigmoid, see, e.g.,~\cite{Aggarwal-18}.

To fully characterise the NN, the weights $w_{j,\ell}^k$ and the biases $b_j^k$ are to be determined via a training procedure in order to learn the form of the regression function from the input data.
This aims to minimise a loss function assessing the mismatch, according to an appropriately defined measure, between the prediction $\Smu[NN][i]$ computed by the net and the \emph{true} data $x_i$, for a set of $\nS$ realisations $(\mu_i,x_i)$.
Common choices of loss functions include the mean absolute error (MAE) and the mean squared error (MSE) defined as
\begin{subequations}\label{eq:loss}
\begin{equation}\label{eq:MAE}
\Lmae := \frac{1}{\nS} \sum_{i=1}^{\nS} \abs{\Smu[NN][i] - x_i} ,
\end{equation}
and
\begin{equation}\label{eq:MSE}
\Lmse := \frac{1}{\nS} \sum_{i=1}^{\nS} (\Smu[NN][i] - x_i)^2 ,
\end{equation}
\end{subequations}
respectively.

The minimisation problem mentioned above is solved using a variety of methodologies inspired by gradient descent including, e.g., stochastic gradient descent and Adam optimiser.
In addition, note that such an optimisation problem commonly requires the introduction of a regularisation term.
As previously mentioned for LS and MLS approximations, response functions based on NNs are not interpolative and $\Smu[NN][i] \neq x_i$, for any $i=1,\ldots,\nS$.
For a detailed introduction to neural networks, interested readers are referred to~\cite{Aggarwal-18}.

%---------------------------------------------------------------------
\section{Proper orthogonal decomposition}
\label{sec:POD}
%---------------------------------------------------------------------

The approximation strategies introduced in the previous section provide the fundamental tools to represent the parametric dependence of the quantities of interest. 
In this section, the main ideas of the proper orthogonal decomposition (POD) are summarised to describe a first technique to construct a reduced basis for the functional approximation of a surrogate model.
Once such a reduced basis is obtained, the approximation techniques discussed in Section~\ref{sec:FuncApprox} can be employed to approximate the dependence of the reduced coefficients on the input parameters.

POD is a dimensionality reduction technique that is employed to construct both purely data-driven surrogate models (POD with interpolation/approximation,  henceforth \emph{PODI}), and physics-based surrogate models (POD with projection, henceforth \emph{POD-RB}) following the framework of reduced basis (RB)~\cite{PateraRozza-07,Quarteroni-QMN-book-15,Rozza-HRS-book-16}.
While both PODI and POD-RB surrogate models rely on reduced spaces constructed from snapshot data, the classification adopted here is mainly algorithmic and focuses on the strategy employed to determine the coefficients of the low-dimensional approximation. In particular, POD-RB methods typically rely on an intrusive projection of the governing equations onto the reduced space, whereas PODI approaches are generally non-intrusive and construct surrogate mappings between the user-defined parameters and the reduced coefficients without explicit access to the operators of the underlying physical problem.
It is worth noting that also hybrid strategies have been proposed in the literature, see, e.g.,~\cite{Rozza-PGR-24}, and that the boundary between these categories may be blurred depending on the specific implementation.

%---------------------------------------------------------------------
\subsection{Dimensionality reduction using POD}
\label{sec:POD-DimRed}
%---------------------------------------------------------------------

Let $\nP=1$ and $\xx_i$,  $i = 1,\dots,\nS$ be a collection of data, obtained either from experiments or high-fidelity simulations,  for $\nS$ values of the input parameter $\mu$.
These vectors, commonly referred to as \emph{snapshots},  are centred and collected in the matrix $\bX$ such that 
\begin{equation}\label{eq:matrixX}
\begin{aligned}
\bX =& \left[ \widehat{\xx}_1\,\,  \widehat{\xx}_2\,\,  \dots \,\, \widehat{\xx}_{\nS} \right],  \\
&\hspace{5pt}\widehat{\xx}_i ={\xx}_i - \bar{\xx}, \text{ and } 
\bar{\xx} =  \dfrac{1}{\nS} \sum_{i=1}^{\nS} {\xx}_{i}.
\end{aligned}
\end{equation}
In some cases, scaling or normalising the snapshots in $\bX$ is recommended to avoid over-representation of some of them, for instance when data features different variables (e.g. velocity and pressure in incompressible flows) with significantly different orders of magnitude.
It is worth noting that the choice of the collection of snapshots in $\bX$ significantly affects the quality of the surrogate models constructed using POD. 
To this end,  Section~\ref{sec:Sampling} reviews some ideas associated with sampling procedures as well as data augmentation strategies to enrich poorly populated datasets.

POD, which is based on Singular Value Decomposition (SVD), is used to eliminate redundant information in data collected in $\bX$. The SVD of $\bX$ reads as
\begin{equation}
\bX = \bUt \bSigma \bVtT,
\end{equation}
where $\bSigma\in\RR^{\ndofx\times\nS}$ is a diagonal matrix containing the singular values of $\bX$ in descending order, $\sigma_1\ge\sigma_2\ge\dots\ge\sigma_\nS\ge0$, whereas $\bUt\in\RR^{\ndofx\times\ndofx}$ and $\bVt\in\RR^{\nS\times\nS}$ are two unitary matrices.

The first $\nS$ columns of $\bUt$ define an orthonormal basis for the linear subspace generated by the snapshots. The most relevant information is contained in the first ones, which correspond to larger singular values. 
In practice,  the first $\nk$ columns of $\bUt$ are kept, $\nk$ being  such that
\begin{equation}\label{eq:criterion}
\sum_{i=1}^\nk \sigma_i \ge (1-\varepsilon) \sum_{i=1}^{\nS} \sigma_i ,
\end{equation}
where the factor $1{-}\varepsilon$ denotes a lower bound for the relative amount of energy captured by the first $\nk$ modes, with $\varepsilon>0$ being a user-defined threshold value.

The resulting reduction is achieved by representing the unknown $\xx$ in terms of only the maintained $\nk$ modes since $\nk$ is expected to be much lower than $\nS$, and $\nS$ is typically taken much lower than the dimension $\ndofx$ of the high-fidelity system. Namely, it holds that
\begin{equation}\label{eq:POD_approx}
\xx \simeq \bar\xx  + \bU \zz =  \bar\xx+ \sum_{i=1}^\nk [\bUt]_i z_i,
\end{equation}
where $\zz=[z_1 \, z_2 \dots z_\nk]^\top$ is the vector of the POD coefficients, and $\bU$ is the matrix containing the first $\nk$ columns of $\bUt$. 

\begin{remark}[High-order SVD]
The strategy described above can be generalised to the case of high-dimensional problems ($\nP > 1$), that is,  when multiple parameters $\bmu$ are involved.
The corresponding reduction approach is based on the high-order SVD~\cite{Lathauwer-LMV-00,Vasilescu-VT-02}, a particular case of orthogonal decomposition proposed by~\cite{Tucker-66}. Whilst this approach provides a reduced description of the high-dimensional space, it is no longer guaranteed that such a representation is unique, with a consequent risk of achieving suboptimal surrogate models.
To circumvent this issue, techniques based on the so-called a posteriori proper generalised decomposition separation are described in Section~\ref{sec:PGD}.
\label{rmrk:HOSVD}
\end{remark}

%---------------------------------------------------------------------
\subsection{PODI}
\label{sec:PODI}
%---------------------------------------------------------------------

Following~\eqref{eq:POD_approx}, for each entry of the dataset $\bX$, the snapshot $\widehat{\xx}_i, \, i=1,\ldots,\nS$ can be expressed in terms of the POD basis as
\begin{equation}\label{eq:PODI}
\widehat{\xx}_i \simeq \bU \zz_i .
\end{equation}
Define the matrix $\bZ = [\zz_1 \,\,  \zz_2\,\,  \dots \,\,  \zz_{\nS} ] \in \RR^{\nk\times\nS}$ collecting the POD coefficients for all $\nS$ snapshots. It follows that
\begin{equation}\label{eq:snapshotApprox}
\bX \simeq \bU \bZ .
\end{equation}
By premultiplying~\eqref{eq:snapshotApprox} by $\bU^\top$, the vector of truncated POD coefficients for all snapshots in the original dataset is obtained, namely,
\begin{equation}\label{eq:PODcoeff}
\bZ = \bU^\top \bX .
\end{equation}

Starting from the reduced-order approximation constructed in~\eqref{eq:POD_approx},  a surrogate model for a new value $\bmuS$ of the parameter is determined by appropriately computing the coefficients in~\eqref{eq:POD_approx}.
PODI, first proposed by~\cite{Ly-LT-01}, employs standard techniques for functional approximation, encompassing both interpolation and regression approaches. 

First, the method expresses $\bmuS$ as a linear combination of the $\nS$ sets of parameters $\bmu \in \mathcal{P}$ used to construct the dataset of high-fidelity solutions. 
This is achieved by interpolating in the parametric space $\mathcal{P}$, that is, by defining the vector $\wwS$ containing $\nS$ weights such that
\begin{equation}\label{eq:paramInterp}
\bmuS = \bTheta \wwS ,
\end{equation}
where each column of the matrix $\bTheta = [\bmu_1 \,\,  \bmu_2\,\,  \dots \,\,  \bmu_{\nS} ] \in \RR^{\nP\times\nS}$ collects one set of parameters associated with the snapshots in $\bX$.
Then, the new set of coefficients $\zzS$ for the reduced-order approximation associated with $\bmuS$ is obtained by interpolating the POD coefficients associated with the $\nS$ snapshots with the weights $\wwS$, leading to
\begin{equation}\label{eq:coeffInterp}
\zzS = \bZ \wwS ,
\end{equation}
and the resulting PODI surrogate model evaluated at point $\bmuS$ follows from~\eqref{eq:POD_approx},  namely,
\begin{equation}\label{eq:PODI surrogate}
\bSmu[POD][\star] := \bar\xx  + \bU \zzS .
\end{equation}

It is worth noting that the above framework can also be extended to the case of non-interpolant functional approximations, such as the regressions approaches reviewed in Sections~\ref{sec:LS},~\ref{sec:MLS}, and~\ref{sec:introNN}.
More precisely, different strategies have been proposed in the literature to perform snapshots interpolation/approximation within PODI, including cubic spline interpolants~\cite{Ly-LT-01},  interpolation to the tangent space to a Grassman manifold~\cite{Amsallem-AF-08}, bivariate interpolation for scattered data~\cite{Alonso-AVV-09}, piecewise linear polynomial functions~\cite{Chinesta-NACC-12},  radial basis functions~\cite{Hassan-WHM-13},  sparse grid collocation~\cite{Pain-XFBPNM-15}, moving least squares~\cite{rama-RSS-16}, least squares on sparse grid points~\cite{Pain-LXFPN-17}, Kriging or Gaussian process regression~\cite{Breitkopf-XBCKSV-10,Hesthaven-GH-18,Pain-XHMFLNGMRP-19}, feed-forward neural networks~\cite{Hesthaven-HU-18},  long short-term memory networks~\cite{Pain-WXFGPG-18},  k-nearest-neighbors regression model~\cite{Willcox-SMPW-19}, and transformer networks~\cite{Pain-WQFFP-22}.
Note that, within the functional approximation framework introduced in Section~\ref{sec:Surrogates}, the different methodologies discussed above should not be interpreted as fundamentally distinct modelling paradigms. 
Indeed, all these approaches represent alternative strategies to address the same underlying problem of approximating the parametric input-output map once a reduced representation of the solution has been introduced. They primarily differ in how the approximation space and the associated surrogate mapping are constructed and exploited and each involves different assumptions, levels of intrusiveness, and trade-off in terms of data requirements, expressiveness, robustness, and generalisation capabilities.
Methods based on geometric interpolation (e.g., interpolation on matrix or Grassmann manifolds) explicitly leverage structural properties of reduced subspaces, whereas regression and learning approaches (e.g., neural networks) approximate the same mappings through generic functional representations.

Of course, when high-dimensional functions are considered,  traditional interpolation and approximation approaches can become computationally unfeasible, requiring an unmanageable amount of data.  
In this context, sparsity-promoting strategies have gained increasing attention stemming from the work~\cite{Brunton-BPK-16}.
Additionally,  in recent years,  research effort has been devoted to data augmentation procedures to enrich poorly populated datasets by devising artificial snapshots introducing relevant information in the training set (see Section~\ref{sec:Augmentation}).

%---------------------------------------------------------------------
\subsection{POD-RB}
\label{sec:PODRB}
%---------------------------------------------------------------------

Assume that data $\xx_i$,  $i = 1,\dots,\nS$ introduced in the previous section are obtained from a known full-order solver (e.g., solving a finite element model or, similarly,  finite volume or finite difference).
Each snapshot is a vector of dimension $\ndofx$ obtained as the solution of the linear system of equations
\begin{equation}\label{eq:fullOrder}
\bK(\bmu) 	\xx(\bmu) = \bff(\bmu),
\end{equation}
where $\bK$ and $\bff$ respectively represent the matrix and the right-hand side vector stemming from the discretisation of the mathematical model under analysis.

An alternative strategy to compute the coefficients $\zz$ in \eqref{eq:POD_approx} for a new value $\bmuS$ of the parameters relies on enforcing that the physical model is verified in the reduced space.
This is achieved via the \emph{projection} of the high-fidelity model \eqref{eq:fullOrder} onto the low-dimensional space spanned by the reduced basis $ \bU$ introduced in \eqref{eq:POD_approx}, that is, by enforcing the physical law of the system in the reduced space.
In this context, the coefficients $\zz$ are selected in order to guarantee, by construction,  compliance of the surrogate model to the physical laws describing the phenomenon under analysis.
Following the RB rationale,  the coefficients $\zz$ in \eqref{eq:POD_approx} are determined via a Galerkin projection of the full-order system \eqref{eq:fullOrder} yielding the POD-RB system
\begin{equation}\label{eq:reducedOrder}
[\bU^\top  \bK(\bmu) \bU ] \, \zz = \bU^\top \bff(\bmu) - \bU^\top  \bK(\bmu) \, \bar\xx ,
\end{equation}
of reduced size $\nk \ll \nS$.

It is worth noting that while POD-RB requires access to the high-fidelity solver employed to perform physical simulations, the resulting surrogate models fulfil the underlying physical laws of the system under analysis.
On the contrary, while PODI allows the construction of surrogate models without communication with the underlying high-fidelity solver, their approximation capabilities are more limited and tend to fail when extrapolation is considered.
This observation motivated an increasing interest towards deep learning solutions, to exploit their generalisation capabilities in the context of POD as mentioned in Section~\ref{sec:PODI}.

If the underlying physics is described by a nonlinear model,  equation \eqref{eq:fullOrder} is to be replaced by a corresponding nonlinear problem
\begin{equation}\label{eq:nonlinFullOrder}
\mathcal{L}_{\bmu} (\xx(\bmu))[ \xx(\bmu) ] = \bff(\bmu),
\end{equation}
where $\mathcal{L}_{\bmu}(\odot)[\odot]$ is the discretisation of a nonlinear PDE (e.g. the Navier-Stokes equations). 
In this context, whilst reduction is possible,  the projection described in \eqref{eq:reducedOrder} requires the evaluation of the nonlinear model \eqref{eq:nonlinFullOrder} using the full-order solution $\xx(\bmu)$.

The introduction of nonlinearity in~\eqref{eq:nonlinFullOrder} naturally highlights several contemporary challenges associated with classical POD-based surrogate models.
In particular, the reliance on linear approximation spaces may lead to slow convergence or loss of accuracy when the solution manifold exhibits pronounced nonlinear features. This limitation can be formally related to an unfavorable decay of the Kolmogorov $n$-width~\cite{Cohen-CD-15}.
In these settings, increasing the dimension of the reduced basis does not necessarily result in accuracy gains, limiting the effectiveness of linear subspace approximations.

These considerations have motivated a growing body of recent work aimed at overcoming linear subspace limitations through nonlinear dimensionality reduction techniques. 
Early examples include kernel principal component analysis (KPCA)~\cite{kPCA-97} and locally linear embedding (LLE)~\cite{LLE-00}, which aim to capture nonlinear structures by embedding the data into low-dimensional nonlinear manifolds. 
Building on these ideas, several approaches have been developed for nonlinear surrogate modeling, such as manifold learning via local maximum-entropy approximants~\cite{Arroyo-MA-13}, manifold walking based on diffuse approximation~\cite{Breitkopf-LRB-15}, and kernel-based extensions of POD~\cite{diez2021nonlinear,Manzoni-SDM-21}. 
While these methods can significantly improve approximation capabilities in the presence of nonlinear solution manifolds, in most existing implementations, they are often intrusive with respect to the high-fidelity model used to generate the snapshots.
This observation has contributed to the growing interest in nonlinear reduced representations based on machine learning. In particular, deep learning approaches exploiting autoencoders aim to construct expressive nonlinear surrogate models while alleviating some of the intrusiveness constraints of classical nonlinear manifold techniques (see Section~\ref{sec:ROMNN}).
In parallel, hybrid approaches combining projection-based methods with data-driven components (e.g., operator learning~\cite{Willcox-PW-16} and closure modelling~\cite{Iliescu-SMLSDI-22}) have been proposed to enhance the robustness and the generalisation capabilities of surrogate models under strong parametric variability.

Overall, these developments reflect a broader shift toward surrogate modelling strategies that seek to balance physical structure~\cite{Ghattas-GW-21}, interpretability and explainability~\cite{Vinuesa-VBM-26} in order to enhance approximation accuracy, computational efficiency, and practical applicability in highly nonlinear and high-dimensional parametric settings.
While a comprehensive review of these emerging directions lies outside the scope of the present work, their discussion here provides essential context for the limitations of linear POD-based models and clarifies the role of nonlinear and data-assisted approaches within the broader functional-approximation framework adopted in this review.

It is worth noting that although nonlinear dimensionality reduction techniques enhance the approximation capability of surrogate models by addressing the limitations of linear subspaces, the need to evaluate the resulting low-dimensional model using the full-order solution $\xx(\bmu)$, see~\eqref{eq:nonlinFullOrder}, significantly limits the computational advantage of projecting the problem onto a reduced space.
In this context, hyper-reduction methods are required to ensure that the resulting nonlinear reduced-order models remain computationally efficient by enabling affordable evaluations of the associated operators.
Available solutions proposed in the literature include, for instance, gappy POD~\cite{Everson-ES-95,Willcox-06}, hyper-reduction collocation~\cite{Ryckelynck-05}, best points interpolation~\cite{Patera-NPP-08}, missing point estimation~\cite{Willcox-AWWB-08}, discrete empirical interpolation method (DEIM)~\cite{Sorensen-CS-10},  least-squares Petrov-Galerkin (LSPG) projection~\cite{Carlberg-CBF-11}, Gauss-Newton with approximated tensors (GNAT)~\cite{Carlberg-CFCA-13}, energy-conserving sampling and weighting (ECSW) method~\cite{Farhat-FACC-14}, and empirical cubature method~\cite{Hernandez-HCF-17}.

Note that classical hyper-reduction approaches are generally \emph{intrusive} since they require access to the high-fidelity solver, thus posing additional constraints in an industrial setting.
To address this issue, different strategies have been proposed.
On the one hand, purely data-driven procedures approximate reduced-order operators from data, without the need to access the operators of the full-order model~\cite{Willcox-PW-16,Willcox-KPW-24}.
On the other hand, deep learning algorithms have been recently exploited for hyper-reduction. For instance, \cite{Hesthaven-CWHZ-21} introduce the physics-reinforced neural network (PRNN) paradigm to circumvent hyper-reduction of nonlinear problems,  \cite{Manzoni-CFM-22} propose to train a neural network replacing DEIM in a nonlinear POD-DEIM algorithm, and a neural empirical interpolation method (NEIM) is discussed in~\cite{Hesthaven-HPH-25}.

From a performance perspective, POD-based surrogate models are typically characterised by very high online efficiency once the reduced basis has been constructed, which makes them particularly attractive for real-time and many-query applications.
However, this advantage is often accompanied by expensive offline stages, requiring intrusive access to the high-fidelity solver in the case of POD-RB, as well as limited expressiveness when the underlying solution manifold exhibits strong nonlinear features or high-dimensional parametric dependence.
As highlighted in the previous discussion, these aspects impose practical limitations on classical POD-based surrogates and motivate the exploration of alternative or complementary strategies.
In particular, learning-based and hybrid data-assisted approaches aim to alleviate some of these constraints by improving flexibility, scalability, and robustness, while still seeking to preserve the efficiency benefits associated with POD-based reduced representations~\cite{Koumoutsakos-WVKS-18,Manzoni-FM-22,Manzoni-BFFM-24,Chinesta-MRGFC-Preprint-24}.

%---------------------------------------------------------------------
\section{Proper generalised decomposition}
\label{sec:PGD}
%---------------------------------------------------------------------

An alternative to POD for the construction of a functional surrogate model relies on the proper generalised decomposition (PGD) framework described in~\cite{Chinesta-Keunings-Leygue}.
Building on the approximation concepts introduced in Section~\ref{sec:FuncApprox}, PGD is here reviewed as a strategy to construct a reduced basis for the surrogate representation of parametric systems. Unlike projection-based approaches such as POD-RB, PGD embeds the parametric dependence directly into a separated representation of the solution. In particular, the idea is to approximate the collection of snapshots by means of a separable representation, that is as the sum of $m$ rank-one terms:
\begin{equation}\label{eq:PGDapprox}
\Xpgd^m := \sum_{j=1}^m \sigma^j \fX^j \otimes \bPsi_1^j \otimes \bPsi_2^j \otimes \cdots \otimes \bPsi_{\nP}^j .
\end{equation}
Each term in \eqref{eq:PGDapprox}, referred to as \emph{mode}, is the product of a vector $\fX^j$ encapsulating the spatial dependence of the data and $\nP$ vectors $\bPsi_k^j$, $k=1,\ldots,\nP$ each depending on one parameter $\mu_k$. Moreover,  $\sigma^j$ denotes the amplitude of the $j$-th mode. 
This framework yields an explicit functional approximation of the quantities of interest while retaining the advantages of the reduced complexity required for efficient surrogate modelling. Similarly to what was discussed for POD, PGD also allows the construction of both physics-based and data-driven surrogate models depending on the strategy employed to compute the modes of the separated representation. In the remainder of this section, different PGD strategies to construct the separated approximation \eqref{eq:PGDapprox} are described.

%---------------------------------------------------------------------
\subsection{A posteriori PGD}
\label{sec:PGDaPosteriori}
%---------------------------------------------------------------------

Assume that the collection of snapshots $\XX$ is constructed according to a tensorial structure, that is, for each parameter $\mu_k$, $k=1,\ldots,\nP$ a matrix $\bX$ (see equation \eqref{eq:matrixX}) is constructed.
The resulting high-dimensional tensor $\XX$ features $\nP{+}1$ dimensions.

The a posteriori PGD, also known as PGD separation in ~\cite{DM-MZH:15} or least-squares PGD in~\cite{PD-DZGH-20}, computes the separated approximation \eqref{eq:PGDapprox} using a greedy approach.
Observing that 
\begin{equation}\label{eq:PGDapproxIncrement}
\Xpgd^m = \sigma^m \fX^m \otimes \bPsi_1^m \otimes \bPsi_2^m \otimes \cdots \otimes \bPsi_{\nP}^m  + \Xpgd^{m-1} , 
\end{equation}
the a posteriori PGD algorithm assumes that $m {-} 1$ modes are known to compute the $m$-th term in the PGD expansion as
\begin{equation}\label{eq:posterioriMin}
\begin{aligned}
&\left(\fX^m , \bPsi_1^m , \ldots , \bPsi_{\nP}^m \right) = \\
&\hspace{15pt}\argmin \left\| \XX {-} \Xpgd^{m-1} {-} \sigma^m \fX^m {\otimes} \bPsi_1^m {\otimes} {\cdots} {\otimes} \bPsi_{\nP}^m \right\|_2 , 
\end{aligned}
\end{equation}
that is, it determines the vectors $(\fX^m , \bPsi_1^m , \ldots , \bPsi_{\nP}^m )$ minimising the residual discrepancy between the full-order data $\XX$ and the previously computed approximation $\Xpgd^{m-1}$ featuring $m{-}1$ terms. 
More precisely,  the greedy procedure aims to compute, at each step, the best approximation $\sigma^m \fX^m \otimes \bPsi_1^m \otimes \cdots \otimes \bPsi_{\nP}^m$ to describe the unresolved part $\XX {-} \Xpgd^{m-1}$ of the target tensor $\XX$. 
From a practical point of view, the nonlinear problem~\eqref{eq:posterioriMin} is solved using an alternating direction scheme and the computation relies on elementary tensorial operations (i.e., products and sums of separated tensors).
A sketch of the a posteriori PGD strategy is reported in algorithm \ref{alg:PGDposteriori}. The method, rigorously presented in~\cite{PD-DZGH-20}, is available in open-source libraries\footnote{\url{https://git.lacan.upc.edu/zlotnik/algebraicPGDtools}\label{ft:PGDmatlab}}\textsuperscript{,}\footnote{\url{https://git.lacan.upc.edu/encapsulated-pgd/demo}\label{ft:ePGD}} providing end-users a black-box allowing the seamless execution of PGD-based operations.
\begin{algorithm}
\caption{The a posteriori PGD algorithm}\label{alg:PGDposteriori}
\begin{algorithmic}[1]
\REQUIRE{For the greedy enrichment loop, the value $\eta^\star$ of the tolerance. For the alternating direction loop, the value $\eta_{\sigma}$ of the tolerance on the amplitude variation and the maximum number of iterations $\niter$.}
\STATE{Compute $\nS$ snapshots.}
\STATE{Set $m \gets 1$ and initialise the amplitude of the spatial mode $\sigma^1 \gets 1$.}

\WHILE{$\sigma^m / \sigma^1 > \eta^\star$}
\STATE{Set $q \gets 1$ and initialise the parametric mode.}
             
\WHILE{$\varepsilon_{\sigma} > \eta_{\sigma}$ or $q < \niter$}

\STATE{Compute the rank-one spatial mode.}

\STATE{Compute sequentially the rank-one parametric modes.}

\STATE{Update the stopping criterion $\varepsilon_{\sigma} = (\sigma^{m,q} - \sigma^{m-1})/\sigma^{m,q}$.}
\STATE{Increase the counter of the alternating direction iterations $q \gets q+1$.}

\ENDWHILE

\STATE{Increase the mode counter $m \gets m+1$.}
\ENDWHILE
\end{algorithmic}
\end{algorithm}

The evaluation of the resulting PGD surrogate model at a new point $\bmuS$ simply requires independently evaluating each parametric mode $\bPsi_k^j$, for $k=1,\ldots,\nP$ and $j=1,\ldots,m$ for the corresponding value of the parameter $[\bmuS]_k$. This leads to a set of scalar coefficients $\mathrm{b}_k^j := \bPsi_k^j([\bmuS]_k)$ acting as weights in the linear combination of the $m$ spatial modes $\fX^j$ as
\begin{equation}\label{eq:PGDsurrogate}
\bSmu[PGD][\star] :=  \sum_{j=1}^m \sigma^j \fX^j \prod_{k=1}^{\nP} \mathrm{b}_k^j .
\end{equation}
Contrary to POD,  this approach is applicable seamlessly for any number $\nP$ of parameters and provides the \emph{optimal} separated representation of the parametric data $\XX$.
In this context, a posteriori PGD represents a natural alternative to the high-order SVD framework discussed in Remark~\ref{rmrk:HOSVD}.
The algorithm employs the amplitude of the modes as a stopping criterion, since it provides a relative measure of the relevance of the newly computed term in the PGD expansion.
Of course, alternative strategies can also be considered, for instance introducing appropriate error control procedures targeting the unknowns of the problem or given quantities of interest, see, e.g.,~\cite{PD-GBCD-17} and~\cite{Chamoin-KCLP-19}.

It is worth noting that the resulting surrogate model~\eqref{eq:PGDsurrogate} is purely data-driven as the governing equations are not enforced in the reduced space. In this context, a posteriori PGD represents an alternative to the PODI solutions discussed in Section~\ref{sec:PODI}, naturally leading to non-intrusive surrogates. Note that, contrary to PODI, this approach does not classify as a \emph{projection-based} reduced order model.

\begin{remark}[Interpretation of the parametric modes]
Note that the parametric modes $ \bPsi_k^j$ are seen both as vectors (e.g. in \eqref{eq:PGDapprox}) and as functions of the variable $[\bmu]_k$ (to compute $\mathrm{b}_k^j$). Frequently, vector  $ \bPsi_k^j$ represents the nodal values of the corresponding function in a 1D grid discretising the $k$-th parametric dimension. In this case, computing $\mathrm{b}_k^j$ consists in a simple interpolation exercise. Nothing prevents, however, adopting a more complex functional representation of the mode, for instance collecting in vector $ \bPsi_k^j$ the coefficients expressing the mode in some functional basis.
\label{remVectFun}
\end{remark}

%---------------------------------------------------------------------
\subsection{Sparsity-promoting PGD}
\label{sec:sparsePGD}
%---------------------------------------------------------------------

In order to reduce the computational cost associated with the construction of the collection of snapshots $\XX$ (especially, in high-dimensional problems),  different strategies to circumvent the need for a tensorial structure of $\XX$ have been discussed in the literature.
The goal is to construct a separated approximation starting from sparse data, not sampling the parametric space according to a Cartesian distribution of snapshots.

A first variant of sparse PGD, known as s-PGD was proposed by~\cite{Ibanez-IAAGCHDC-18}. This strategy relies on constructing a PGD separated expansion approximating the norm in~\eqref{eq:posterioriMin}  using a collocation approach.
In this context, the limited amount of data combined with high-order of polynomial approximations was responsible for overfitting phenomena. To avoid this issue,  s-PGD employs a degree adaptive strategy, forcing the degree of the first PGD modes to remain low and progressively increasing it, when higher order variations are sought by the PGD solution.

Following a similar rationale,  \cite{Gravouil-LBG-18} proposed a least-squares PGD strategy based on adaptive sparse grids~\cite{Brumm-BS-17,Piazzola-PT-24}. The idea is to construct a PGD separated approximations minimising only the error on sampling points selected according to sparse grids, thus allowing for parsimonious surrogate models of high-dimensional parametric problems.
Note that, contrary to~\cite{Ibanez-IAAGCHDC-18}, the PGD modes computed by \cite{Gravouil-LBG-18} are determined exclusively from the minimisation of the residual in~\eqref{eq:posterioriMin}, without any additional constraint.

In~\cite{Rocas-RGZLD-21},  the standard (Euclidean or Frobenius) norm in~\eqref{eq:posterioriMin} is replaced by a norm in the parametric space $\mathcal{P}$. Moreover,  the sampling points selected to populate the dataset are used as integration points, with proper discrete weights, to compute the newly introduced norm. Thus, the weights of the quadrature are adapted to optimise the integration order in a multidimensional Newton-Cotes fashion.
It is worth noting that the functions selected in this work for the spatial and parametric approximation of~\eqref{eq:PGDapprox} were standard $\mathcal{C}^0$ finite element shape functions.
Whilst these functions are more stable than high-order polynomials and yield a sparser system owing to their compact support and small stencil,  they might suffer from a lack of smoothness (e.g., jumps in first-order derivatives and singularities in second-order derivatives).
Hence, the sparse PGD introduces a regularisation term during the computation of the separated approximation~\eqref{eq:PGDapprox} in order to penalise the Euclidean norm of the gradient of the approximation.

Finally,  \cite{Chinesta-SCCC-23} discuss two variants of the s-PGD, namely the regularised sparse PGD (rs-PGD), and the doubly sparse PGD (s2-PGD).
On the one hand, the rs-PGD combines Ridge and LASSO regressions by introducing two regularisation terms, one based on the $\eltwo$ norm and one based on the $\elone$ norm of the unknowns: the former provides a smoothing of the PGD response surface to avoid overfitting, while the latter sparsifies the solution retaining only the contributions of the most relevant terms in the approximation.
On the other hand, the s2-PGD promotes sparsity by enforcing LASSO regularisation in each PGD direction independently.

%---------------------------------------------------------------------
\subsection{A priori PGD}
\label{sec:PGDaPriori}
%---------------------------------------------------------------------

Similarly to the rationale discussed in Section~\ref{sec:PODRB} for POD, also PGD-based surrogate models can be constructed enforcing compliance of physical laws at the reduced level.
This framework leads to the so-called \emph{a priori} PGD rationale, that is, the construction of the rank-one approximations $\fX^m , \bPsi_1^m , \ldots , \bPsi_{\nP}^m$ in the separated representation~\eqref{eq:PGDapprox} using a greedy algorithm that sequentially minimises the residual of the underlying governing equation.
For a linear equation, this yields, for all $\bmu \in \mathcal{P}$, the solution of the problem
\begin{equation}\label{eq:PGDlinSys}
\Kpgd(\bmu) \Xpgd^m = \Fpgd(\bmu) ,
\end{equation}
where $\Kpgd$ and $\Fpgd$ are the so-called \emph{separated tensors}, see~\cite{PD-DZGH-20}, defined as
\begin{subequations}\label{eq:PGDdata}
\begin{equation}\label{eq:PGDdataK}
\Kpgd(\bmu) :=  \sum_{j=1}^{M_K} \bK^j \xi^j(\bmu) = \sum_{j=1}^{M_K} \bK^j \prod_{k=1}^{\nP} \xi_k^j(\mu_k),
\end{equation}
\begin{equation}\label{eq:PGDdataB}
\Fpgd(\bmu) :=  \sum_{j=1}^{M_f} \bff^j \eta^j(\bmu) = \sum_{j=1}^{M_f} \bff^j \prod_{k=1}^{\nP} \eta_k^j(\mu_k) ,
\end{equation}
\end{subequations}
with $\bK^j$ and $\xi^j$ (respectively, $\bff^j$ and $\eta^j$) being the $j$-th spatial and parametric modes of the separated representation associated with the matrix $\bK(\bmu)$ (respectively, the vector $\bff(\bmu)$) in~\eqref{eq:fullOrder}.

Contrary to a posteriori PGD, this approach does not rely on the availability of precomputed snapshots.
The numerical strategy for a priori PGD, sketched in algorithm~\ref{alg:PGDpriori}, is based on the alternating direction method to minimise the residual of the parametric problem~\eqref{eq:PGDlinSys}.
The user is not required to select any snapshot in order to construct the reduced model: the greedy algorithm is in charge of computing the full-order solutions to minimise the residual of the system during the enrichment process, and the required number $m$ of terms is automatically determined according to a user-defined tolerance.
This method has been successfully applied to many (parametric) physical systems over the years and interested readers are referred to the review article~\cite{Chinesta-CLC-11}.
Hence, whilst a posteriori PGD highly depends on the initial dataset of solutions in order to construct an accurate surrogate model, a priori PGD circumvents this issue since it does not require any previous knowledge of the solution.
Note that, given the greedy nature of the enrichment procedure, a priori PGD computes modes sequentially.
On the contrary, similarly to POD, the a posteriori PGD procedure can be easily parallelised as snapshots computation is independent of one another, significantly reducing the computational cost of the offline phase.
Moreover, as previously noticed, the purely data-driven framework of a posteriori PGD allows to devise non-intrusive implementations independent of the underlying physical problem.
On the contrary, a priori PGD requires knowledge of the equations governing the physical system, that is, access to the matrices $\bK^j, \ j=1,\ldots,M_K$ and vectors $\bff^j, \ j=1,\ldots,M_f$ associated with the spatial discretisation of the phenomena under analysis.
In this context, most of the early implementations of PGD were intrusive with respect to the method employed for the discretisation of the spatial problem, requiring full access to the high-fidelity solver and yielding a significant limitation to the applicability of the method in realistic, industrial settings, where the use of commercial software is often preferred. 
More recently,  significant effort has been devoted to develop non-intrusive PGD-ROMs integrated with commercial and open-source libraries such as SAMCEF~\cite{Ladeveze-CNLB-16}, Abaqus~\cite{zou2018nonintrusive}, OpenFOAM~\cite{tsiolakis2020nonintrusive,Tsiolakis-TGSOH-22}, VPS/Pamcrash~\cite{Rocas-RGZLD-21}, and MSC-Nastran~\cite{Cavaliere-CZSLD-21,Cavaliere-CZSLD-22}.

\begin{algorithm}
\caption{The a priori PGD algorithm}\label{alg:PGDpriori}
\begin{algorithmic}[1]
\REQUIRE{For the greedy enrichment loop, the value $\eta^\star$ of the tolerance. For the alternating direction loop, the number of iterations $\niter$.}
\STATE{Set $m \gets 1$ and initialise the amplitude of the spatial mode $\sigma^1 \gets 1$.}

\WHILE{$\sigma^m / \sigma^1 > \eta^\star$}
\STATE{Set $q \gets 1$ and initialise the parametric prediction.}
             
\WHILE{$q < \niter$}

\STATE{Contract the spatial direction.}
\STATE{Project in the parametric direction and assembly the parametric system (sequentially, for all parameters).}
\STATE{Solve the linear system for the parametric mode.}

\STATE{Contract the parametric direction.}
\STATE{Project in the spatial direction and assembly the spatial system. }
\STATE{Solve the linear system for the spatial mode.}
             
\STATE{Increase the counter of the alternating direction iterations $q \gets q+1$.}

\ENDWHILE

\STATE{Increase the mode counter $m \gets m+1$.}
\ENDWHILE
\end{algorithmic}
\end{algorithm}

In order to provide a fully non-intrusive framework completely independent of the underlying physical model, \cite{PD-DZGH-20} propose to devise the PGD surrogate model as a low-rank approximation of a high-dimensional discretised equation~\eqref{eq:fullOrder}.
The resulting a priori PGD method is available in open-source libraries\textsuperscript{\ref{ft:PGDmatlab},\ref{ft:ePGD}} which provide black-box tools in MATLAB, Python, and C++ to perform operations on separated tensors and seamlessly construct PGD-based surrogate models.
For a detailed critical and computational comparison of a priori and a posteriori PGD methodologies, interested readers are referred to~\cite{MG-GBSH-21}.

To summarise, PGD-based approaches construct parametric surrogates with explicit functional dependence on the parameters, which can offer compact representations and extremely efficient online performance.
Nevertheless, their performance is closely tied to the separability assumption of the underlying operators.
Indeed, while separated representations can significantly reduce computational complexity, their effectiveness relies on the existence of low-rank, globally separable structures, which may be difficult to achieve in the presence of strong nonlinearities or tightly coupled parameters. 
These considerations naturally motivate the exploration of more flexible surrogate constructions, including learning-based approaches.
In particular, recent developments in this field include a PGD algorithm coupled with the shooting method to treat nonlinear dynamical systems~\cite{Park-LLP-24}, as well as hybrid approaches exploiting deep learning to enhance PGD-based surrogate models. These include a strategy to learn tensor decomposition~\cite{Liu-LMGLL-24}, a hybrid PGD-NN algorithm combining neural interpolation with tensor-decomposed surrogate structures~\cite{Genet-DSG-25}, and a learning-assisted technique to construct a PGD-inspired separated representation by learning the modes of the solution via the training of NNs~\cite{Pardo-BTUP-25}.

%---------------------------------------------------------------------
\section{Neural networks}
\label{sec:NN}
%---------------------------------------------------------------------

The methodologies discussed in the previous sections rely on building structured reduced representations of the solution manifold. In contrast, this section presents an overview of recent approaches based on neural networks to approximate the relation between input parameters and output quantities of interest.
This introduces a complementary paradigm in which surrogate models are constructed as \emph{learned} functional representations.
In this context, NNs can be employed in two different ways.
On the one hand, they can be used within the functional approximation framework introduced in Section~\ref{sec:FuncApprox}, learning the parametric dependence of quantities of interest starting from reduced representations such as those obtained via POD or PGD.
On the other hand, NNs can be employed to directly learn nonlinear reduced representations of the data itself, thus providing alternative dimensionality reduction mechanisms.
These approaches can be interpreted, in many cases, as methodological responses to overcome the limitations identified in the projection-based and separated-representation techniques discussed in Section~\ref{sec:POD} and Section~\ref{sec:PGD}, respectively. These include the previously reviewed shortcomings due to linear subspace restrictions, algorithmic intrusiveness, and scalability in nonlinear and high-dimensional parametric regimes.

The main idea explored in this section is to approximate the collection of data (i.e. snapshots) using a generic, possibly \emph{black-box}, functional representation based on the NN paradigm introduced in section~\ref{sec:introNN}.
Within this perspective, NNs act as universal approximators capable of capturing complex nonlinear dependencies and complement the surrogate modelling strategies introduced earlier.
As previously done for POD and PGD strategies, the fundamental ideas to introduce physical knowledge in the resulting surrogate models are also explored, following the notions of learning bias and inductive bias.

%---------------------------------------------------------------------
\subsection{Dimensionality reduction using neural networks}
%---------------------------------------------------------------------

Neural networks performing dimensionality reduction can be constructed by setting the output data $\Sae$ to represent an approximation of the \emph{true} data $\xx$ provided as input.
In this context, a low-rank approximation is obtained by compressing the information in a lower-dimensional manifold known as latent space.
Note that this is done in an unsupervised context, where no piece of information describing dependence on parameters $\bmu$ is available, only snapshots $\xx_i$ for $i=1,\ldots, \nS$ which are seen both as input and output.

\begin{figure*} [h]
	\centering	
	\includegraphics[width=0.9\textwidth]{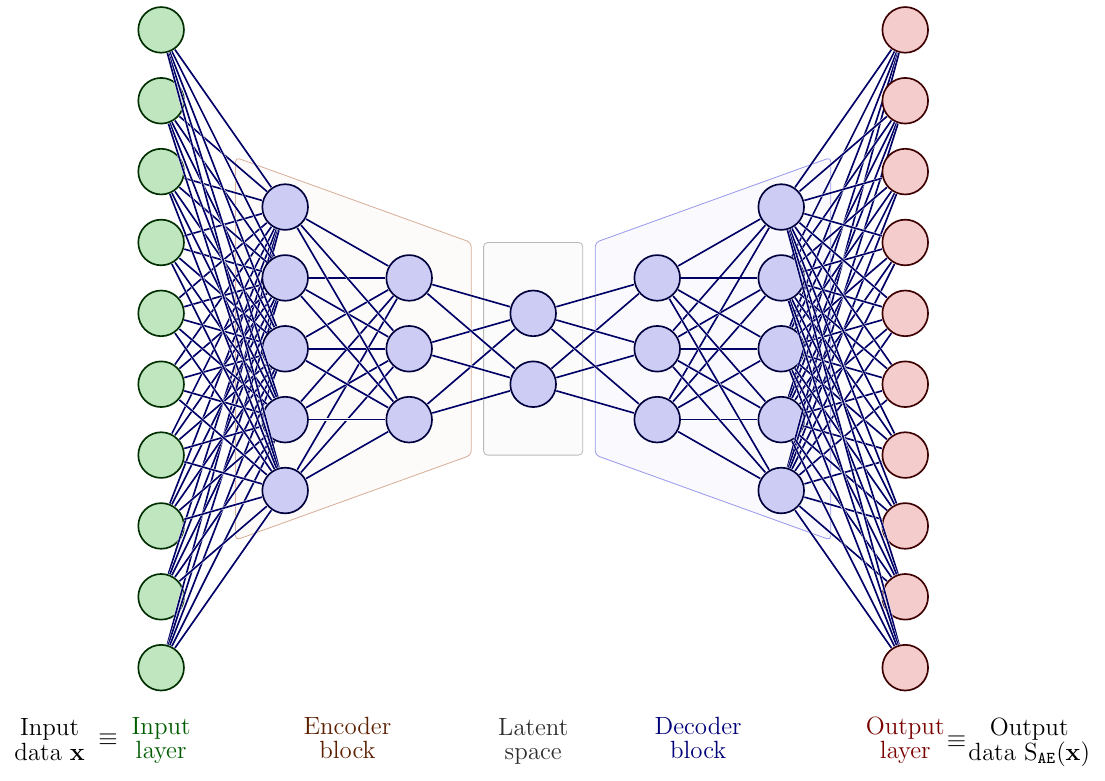}
	\caption{Illustration of an autoencoder neural network.}
	\label{fig:AE}
\end{figure*}

The architecture (see Figure~\ref{fig:AE}) of the resulting net, known as autoencoder,  features a compression block (viz., the \emph{encoder}) reducing the dimensionality of the input data to a low-dimensional \emph{latent space} and a decompression block (viz., the \emph{decoder}) upscaling the information to a space of the same dimension as the input. 

The training of an autoencoder aims to determine the weights and the biases of the net by minimising the mismatch between the reconstructed output $\Sae[i]$ and the input $\xx_i$ for $i=1,\ldots, \nS$, for instance, measured using the MAE and MSE loss functions normalised with respect to the input size $\ndofx$
\begin{subequations}\label{eq:lossAE}
\begin{equation}\label{eq:aeMAE}
\Lmae^{\texttt{AE}} := \frac{1}{\nS \, \ndofx} \sum_{i=1}^{\nS} \abs{\Sae[i] - \xx_i} ,
\end{equation}
and
\begin{equation}\label{eq:aeMSE}
\Lmse^{\texttt{AE}} := \frac{1}{\nS \, \ndofx} \sum_{i=1}^{\nS} \|\Sae[i] - \xx_i\|^2 .
\end{equation}
\end{subequations}

In~\cite{Hinton-AE-06},  the potential superiority of autoencoders with respect to principal component analysis (PCA) for dimensionality reduction has been numerically investigated, whereas interested readers are referred to~\cite{DimReductionReview-21} for a detailed comparison of linear and nonlinear dimensionality reduction algorithms.

%---------------------------------------------------------------------
\subsection{Data-driven and physics-based approximations}
\label{sec:funcNN}
%---------------------------------------------------------------------

This section reviews deep learning contributions to perform functional approximation in the context of computational science and engineering applications following purely data-driven approaches as well as accounting for physical knowledge.

In~\cite{Sevilla-BSHM-21}, a data-driven approach to NN training is discussed. The authors propose a multi-output feed-forward NN trained to predict the response surface of quantities of engineering interest depending on the problem parameters. The results highlight the superior performance of the method with respect to POD, in the presence of solutions with low regularity (e.g., shocks).

Nonetheless, it is known that data-driven models might fail to provide physically-admissible solutions, especially when extrapolating outside of the training set.
Similarly to POD-RB and a priori PGD, physical information can be encapsulated also in machine learning-based surrogate models, especially using NNs.
To the best of the authors' knowledge, existing works on functional approximation based on deep learning can be classified in two groups, according to the strategy employed to include physical knowledge in the NN framework, namely \emph{learning bias} and \emph{inductive bias}.

Methods introducing learning bias incorporate physical knowledge during training, see, e.g., the DeepONet paradigm in~\cite{Karniadakis-LJPZK-21} and~\cite{Perdikaris-WWP-21}, where (parametric) nonlinear operators are learned exploiting the universal approximation theorem.
A similar approach is also discussed in~\cite{Anandkumar-LKALBSA-20}, with the so-called Fourier neural operators (FNO).
To reduce the dimensionality of the resulting problem,  \cite{Azizza-RRA-22} couples a U-Net architecture~\cite{UNet-15} to reduce the dimensionality with a neural operator, whereas \cite{Karniadakis-KGKS-24} proposes to train the neural operator on the latent space, by coupling DeepONets with autoencoders.
Other examples of hybrid modelling blending physics and data via learning bias are presented in~\cite{Farhat-AF-23,VanDerMeer-RKV-23}.

Strategies relying on inductive bias introduce physical information in the NN architecture, for instance enforcing fulfilment of thermodynamics laws.
In~\cite{Cueto-HBGCC-21-JCP}, a feed-forward net is enriched with the GENERIC framework~\cite{GENERIC-97} yielding a NN that preserves the structure of the underlying physical system.
Similarly to what was discussed for DeepONets, the structure-preserving NN mentioned above can be coupled with a sparse autoencoder to reduce the dimensionality of the problem and achieve a parsimonious representation of the system~\cite{Cueto-HBGCC-21-CMAME}.

%---------------------------------------------------------------------
\subsection{Surrogate models based on neural networks}
\label{sec:ROMNN}
%---------------------------------------------------------------------

This section focuses on the use of NNs for the dimensionality reduction required to construct reduced-order and surrogate models.
Note that this goes beyond the discussion of Section~\ref{sec:POD} where the application of NNs was mainly targeting hyper-reduction of nonlinear problems and interpolation in the reduced space generated by classical projection ROMs.

The first proposals in this direction can be traced back to~\cite{Carlberg-LC-20}, where an autoencoder using convolutional layers is employed for dimensionality reduction, yielding nonlinear surrogates (namely, manifold Galerkin ROM and manifold least-squares Petrov Galerkin ROM) inspired by classical Galerkin and LSPG projection ROMs and to~\cite{Duraisamy-XD-20}, linking the autoencoder to a feed-forward net mapping the problem parameters to the latent space representation.
A similar strategy is also discussed in~\cite{Manzoni-FDM-21}.
Stemming from these works, several aspects related to nonlinear dimensionality reduction for nonlinear problems were explored, including hyper-reduction techniques for shallow masked autoencoders~\cite{Choi-KCWZ-22} and for deep convolutional autoencoders~\cite{Rozza-RSR-23}.

In addition, different strategies have been proposed to perform functional approximation upon reducing the dimensionality of the problem using convolutional autoencoders.
In~\cite{Maulik-MLB-21}, the authors employ autoencoders to compress the information in the latent space and long short-term memory (LSTM) networks~\cite{LSTM-97} for interpolation, as commonly adopted in data-driven surrogate modelling of nonlinear dynamical systems (see, e.g.,~\cite{Vlachas-VBWSK-18}), whereas~\cite{Pain-WGPQFP-21} extends this model by maintaining LSTMs and introducing a self-attention module in the convolutional autoencoder to enhance the ability of the network to extract nonlocal features.
In~\cite{Bouklas-KBCOYB-22}, autoencoders are combined with radial basis functions,  \cite{Farimani-HB-23} proposes to learn the dynamics in the reduced encoded space using a transformer network~\cite{Transformers-17},  and \cite{Wang-WZZ-24} uses a convolutional autoencoder for compression and a convolutional NN to map the problem parameters to the latent space.
In~\cite{MG-GH:26}, the authors combine feed-forward NNs with a decoder block to learn the map from the problem parameters to the solution of an outer-loop application in parametrised topology optimisation of elastic structures.
An alternative approach to leverage simulation data defined on unstructured grids is proposed in~\cite{Bouklas-KBOCBY-22}, where convolutional layers are replaced by a feed-forward NN combined with a self-supervised learning strategy of Barlow Twins~\cite{BarlowTwins-21} to compress information.

More recently, other NN architectures have been considered, such as variational autoencoders combined with projection-based POD~\cite{Chatzi-SVGDC-24} and $\beta$-variational autoencoders~\cite{VAE-13} to extract low-dimensional representations of the system that can be coupled with convolutional NNs~\cite{Vinuesa-ELHV-22} and transformers~\cite{Vinuesa-WSVV-24}. 
Moreover, generative adversarial nets~\cite{GAN-14} have been proposed to devise generative ROMs exploiting autoencoder-based dimensionality reduction~\cite{Rozza-CDR-24}, whereas probabilistic surrogates for forecasting the dynamics of high-dimensional complex systems were presented in~\cite{Koumoutsakos-GKK-24}.

Recent developments in NN-based surrogate modelling increasingly focus on key challenges related to generalisation, robustness, and scalability in high-dimensional parametric settings.
Indeed, although autoencoders have shown significant generalisation capabilities, interpretability of the latent space often remains a challenge. This is particularly critical in the context of the applicability of such techniques to realistic industrial scenarios as understanding and controlling the parameters involved in a complex system is of utmost importance. 
For a discussion on this topic and solutions to achieve interpretable latent spaces, interested readers are referred to~\cite{Kutz-CLKB-19,Carlberg-LC-20,Choi-FHC-22}.
Similarly, explainability has become an increasingly important objective, with growing attention devoted to feature-attribution and sensitivity-analysis techniques aimed at understanding how learned representations and input parameters influence model predictions, thereby improving transparency and trust in surrogates constructed using machine learning techniques~\cite{Lundberg-LL-17,Vinuesa-CHDQLLMHLMV-24}.

%---------------------------------------------------------------------
\subsection{Surrogate models based on graph neural networks}
\label{sec:ROMGNN}
%---------------------------------------------------------------------

Standard NN approaches applied to surrogate models do not exploit the spatial correlation of data naturally encapsulated in the underlying discretisations used to compute the snapshots of the full-order problem.
To address this limitation and allow, for instance, employing data from unstructured and deforming grids, the paradigm of graph neural networks (GNNs) proposed in~\cite{GNN-08} has gained increasing attention in recent years.

Starting from the work on \texttt{MeshGraphNets}~\cite{Battaglia-PFSB-20,Manzoni-FFTM-23}, several studies have explored the suitability of GNNs to simulate physical problems while accounting for the underlying geometric configuration of the domain under analysis.

In the context of data-driven surrogate models, \cite{Gunzburger-GGJW-22} showed that graph autoencoders outperform standard autoencoders, allowing the use of unstructured data. It is worth highlighting that, from the presented results, authors conclude that standard POD-based ROMs still offer competitive solutions when the latent space is sufficiently large.

Recent works also explored the suitability of enriching GNNs with physical knowledge of the underlying system. 
In~\cite{Cueto-HBCC-22}, the authors propose to include fulfilment of thermodynamics laws in the GNN architecture via the GENERIC framework~\cite{GENERIC-97}.
Similar work has been presented in~\cite{Cueto-BBGC-25} exploiting the \texttt{MeshGraphNets} architecture.
In~\cite{Pichi-PMH-24}, a surrogate model is devised by coupling a graph-based encoder using convolutional layers, a multilayer perceptron learning the map between the physical parameters and the latent space and a decoder to retrieve the high-dimensional information.
The generalisation capabilities of GNNs for different geometries were also exploited in~\cite{Neron-MVFN-24} to train surrogate models learning a separated PGD structure for geometrically parametrised systems.

\begin{remark}[Data availability]
A crucial aspect for the construction of surrogate models based on POD, PGD, or NNs as discussed in Sections~\ref{sec:POD},~\ref{sec:PGD}, and~\ref{sec:NN}, respectively, is the availability of a \emph{sufficient amount} of high-quality data.
In realistic, many-query scenarios, high-fidelity data are generally scarce and need to be complemented with information of alternative sources, possibly noisy and less reliable.
In this context, optimising the choice of high-fidelity data and devising robust techniques capable of exploiting information from heterogeneous sources and with different levels of credibility is of paramount importance.
These considerations naturally motivate the review of the following sections on multi-fidelity modelling and adaptive dataset sampling strategies as key enabling technologies to enhance the efficiency, robustness, and scalability of surrogate models in complex parametric settings.
From a performance standpoint, these techniques play a central role in surrogate modelling, as they directly govern training cost, data efficiency, and the ability to meet accuracy requirements under limited computational budgets, independently of the specific surrogate construction methodology.
\end{remark}

%---------------------------------------------------------------------
\section{Multi-fidelity models}
\label{sec:MultiFidelity}
%---------------------------------------------------------------------

The techniques described in the previous sections starts from the assumption of having high-fidelity data available for the construction of surrogate models. Nonetheless, in engineering practice, it is common to have access to multiple sources of data, both experimental and simulation-based, with different levels of credibility. In general, data include few, reliable high-fidelity results and many, low-fidelity, possibly noisy, ones. In this context, techniques to concurrently leverage different data types to construct \emph{multi-fidelity} surrogate models have gained increasing attention in recent years, see the review papers~\cite{Park-PHK-17-reviewMF,Gunzburger-PWG-18,Breitkopf-KBBD-22-reviewMF}.

The idea of multi-fidelity methods is to exploit multiple datasets of different fidelities to construct a surrogate model that better represents the phenomenon under analysis than a reduced model built from a single dataset, while reducing the reliance on (and the computational cost of) high-fidelity data alone.
In a computational framework, the datasets mentioned above are obtained via simulations of different fidelities, with high-fidelities (HiFi) featuring accurate but costly and rarer solutions, and low-fidelities (LoFi) providing inexpensive evaluations which might lack of sufficient resolution to capture the features of the problem.
Multi-fidelity methods aim at leveraging LoFi data to extract the major trends of the phenomenon under analysis and exploit HiFi data to locally refine the solution when increased accuracy is required. The correlation between the employed datasets allows to infer trends from LoFi data and exploit it to approximate the HiFi surrogate model, especially for configurations in which HiFi data are scarce.

%---------------------------------------------------------------------
\subsection{Multi-fidelity models based on interpolation}
\label{sec:MultiInterp}
%---------------------------------------------------------------------

Stemming from the work by Kennedy and O'Hagan on co-Kriging~\cite{CoKriging-00}, several authors have explored strategies to capture the discrepancy among models of different fidelities using Gaussian processes.
Co-Kriging has been further enriched with a variance estimator, see~\cite{Keane-FSK-07}, to guarantee robustness to varying degrees of noise on the different levels of fidelity considered.
In~\cite{Gortz-HG-12}, a hierarchical model is presented by leveraging a Kriging interpolation of LoFi data to describe the trend of the quantity of interest, while the variation in the LoFi space is mapped to the HiFi data to improve the surrogate representation.
Moreover,  information on the gradient of the response surface can be integrated in the construction of the Kriging model via the gradient-enhanced Kriging approach discussed in~\cite{Zimmermann-HGZ-13}.

Another class of multi-fidelity approaches leverages the multigrid rationale of employing a hierarchy of grids with different resolutions to accelerate convergence, appropriately transferring information between coarse and fine grids,  to define a multi-fidelity framework for high-dimensional problems.
The first proposal in this direction is the multilevel Monte Carlo (MLMC) method~\cite{Giles-08}. MLMC introduces multiple levels of approximation (many low-cost LoFi data and few expensive HiFi data) and exploits a telescoping sum identity to express the expectation of the HiFi simulation as the sum of the differences between successive levels of approximation of lower fidelity.
In~\cite{Willcox-NW-14},  the authors devise a multi-fidelity strategy in which the estimator from the low fidelities is used to reduce the sampling variance of the estimator from the high fidelity, leading to a method closely related to MLMC, with suitably defined control variate estimators~\cite{Nelson-87}. 
Following the same rationale of a hierarchy of spatial approximations,  the multilevel stochastic collocation (MLSC) and the multilevel quasi-Monte Carlo (MLQCM) were introduced in~\cite{Gunzburger-TJWG-15} by employing quadrature formulae based on sparse grids and quasi-Monte Carlo on each level.
An alternative approach to multi-fidelity via stochastic collocation is discussed in~\cite{Xiu-NGX-14}, where the authors first sample a large number of LoFi configurations and then employ a greedy procedure to select a limited number of points to evaluate the HiFi model and use them as collocation points for the multi-fidelity surrogate.
Similarly,  the multi-index Monte Carlo (MIMC)~\cite{Nobile-HNT-16} and the multi-index stochastic collocation (MISC)~\cite{Nobile-HNTT-16} respectively employ Monte Carlo and sparse grids quadratures and apply the telescopic sum identity mentioned above to the multi-index of the discretisation. The introduction of the multi-index allows to significantly reduce the variance of the hierarchical differences, while achieving a discretisation for each parameter independent of the others.
Finally,  \cite{Piazzola-PTPBSD-23} presents a benchmark test of a complex naval engineering setting to compare accuracy and performance of surrogate models based on MISC and stochastic radial basis functions~\cite{Diez-VDGSICCS-15}, where RBFs are characterised by a stochastic kernel to be determined during a training step.
An alternative multi-fidelity strategy based on deterministic RBFs, see~\cite{Lv-SLSZ-19}, consists of adjusting the predictions of a RBF-based model using LoFi data to match the HiFi outputs via suitably-defined correlation matrix and scaling factors.

\begin{remark}[Multi-fidelity surrogates with interpolation]\label{rmrk:ROM-LoFI-1}
Whilst the datasets of different fidelities can be obtained according to several strategies (e.g., using different levels of mesh refinement in the full-order computation), ROMs offer a natural framework to construct LoFi data of parametric systems.
In this context, several works have explored the coupling of dimensionality reduction techniques with the strategies described above towards the definition of multi-fidelity ROMs.
In~\cite{Mifsud-MMS-16}, POD is combined with the variable-fidelity model to construct a reduced basis by performing SVD on a dataset combining LoFi and HiFi snapshots selected at fixed sample points in the domain.
An alternative approach aims to compute a bridging function to map the POD coefficients obtained from LoFi snapshots to the space of mode coefficients associated with HiFi data~\cite{Wang-WKZ-20}.
Finally, POD has been combined with Kriging~\cite{Breitkopf-BBLSV-17,Zimmermann-BOZ-18} and co-Kriging~\cite{Breitkopf-XZBVZ-18} to construct response surfaces in a space of reduced dimensionality using data of different fidelities.
\end{remark}

%---------------------------------------------------------------------
\subsection{Multi-fidelity models based on regression}
\label{sec:MultiRegr}
%---------------------------------------------------------------------

Non-interpolatory multi-fidelity models exploit the rationale of regression-type approaches described in Section~\ref{sec:FuncApprox}.
In~\cite{Zhang-ZKPH-18}, linear regression is performed using the LoFi data as basis function and minimising the discrepancy between HiFi and scaled LoFi data by means of a least-squares regression.
This framework has been subsequently extended to employ a moving least-squares approach~\cite{Wang-WLZYLS-21} to simultaneously compute the LoFi scaling factors and the coefficients of the discrepancy function to ensure accurate data fusion.

Regression-type multi-fidelity models have also been enhanced by means of Bayesian inference. 
In~\cite{Yan-YZ-19}, a multi-fidelity model based on polynomial chaos expansion is constructed, with LoFi data employed for an initial broad exploration and HiFi data used in the critical regions where the posterior approximation is insufficient.
Similarly,  Bayesian inference has been incorporated in the estimation of regression parameters~\cite{Bessa-YCB-Preprint-24}.  In this context,  the authors propose to combine a deterministic model for regression using LoFi data (e.g.,  kernel ridge regression) with a Bayesian model based on Gaussian process regression (GPR) to account for the discrepancy between LoFi and HiFi data. This work also explores a setup combining deep and Bayesian NNs for the deterministic and statistical regression, respectively.
Closely related to kernel ridge regression, support vector regression~\cite{Vapnik-VGM-96} and variations of such a technique have been successfully employed to identify and describe the relation between LoFi and HiFi data, after introducing a kernel function to map the input space into a high-dimensional feature space~\cite{Lv-SLSS-20,Lv-SLX-23}.

In order to provide multi-fidelity predictions with uncertainty estimates, Gaussian processes have been commonly employed in the literature.
Proposed strategies differ in the modelling of the correlation between the fidelity levels: linear models of coregionalisation assume that LoFi and HiFi data are linear combinations of shared latent Gaussian processes~\cite{Garland-GLRD-20}; nonlinear models of coregionalisation, see, e.g., ~\cite{Kirby-WXKZ-21}, perform nonlinear transformations to propagate the latent basis functions across fidelity levels and employ a matrix Gaussian process prior for the corresponding weights, going beyond the linear assumption; deep Gaussian processes~\cite{Hebbal-HBTM-21} rely on a hierarchy of layers of Gaussian processes mapping LoFi to HiFi data by means of nonlinear transformations; the nonlinear autoregressive Gaussian process regression (NARGP)~\cite{Karniadakis-PRDLK-17} uses a Gaussian process to model the discrepancy between fidelities, learning the form of the associated nonlinear function from data.
For a recent review on the subject, interested readers are referred to~\cite{Do-DZ-23}.

\begin{remark}[Multi-fidelity surrogates with regression]
As previously mentioned in Remark~\ref{rmrk:ROM-LoFI-1}, ROMs can be employed to define LoFi data within multi-fidelity settings.
In~\cite{Guo-KGH-20}, the authors employ LoFi data to identify the parameter locations of the required HiFi simulations by means of a nonlinear autoregression Gaussian process and construct a POD basis using the HiFi snapshots.
To reduce the dimensionality of the parametric space,  \cite{Rozza-RTMOR-23} proposes to fuse LoFi models, obtained by means of active subspaces~\cite{Constantine-15} or nonlinear level-set learning~\cite{Zhang-ZZH-19},  with HiFi data exploiting the NARGP paradigm described above.
The authors further build on top of this idea to combine active subspaces to reduce the dimensionality of the parametric space,  NARGP to model the relation between fidelities, with LoFi data guiding the selection of HiFi computations,  and PODI with RBF to reduce the dimensionality of the solution space constructed from HiFi data~\cite{Rozza-TFSSR-23}.
\end{remark}

%---------------------------------------------------------------------
\subsection{Multi-fidelity models based on neural networks}
\label{sec:MultiNN}
%---------------------------------------------------------------------

Recently, there has been a growing interest towards the use of NNs to devise multi-fidelity surrogates.
In~\cite{Karniadakis-MK-20}, a deep NN architecture is devised by combining a net trained on LoFi data to learn a coarse approximation and two nets aiming to capture linear and nonlinear correlations between LoFi and HiFi data. This framework is further extended in~\cite{Karniadakis-MBK-21} to account for uncertainty quantification via Bayesian inference techniques and in~\cite{Karniadakis-HPKS-23} to learn operators using the DeepONet framework.
Similarly,  a three nets architecture based on convolutional NNs is discussed in~\cite{Zhang-ZXJZZ-21} to integrate data of two fidelities, dynamically incorporating  in the dataset new HiFi data points based on previously computed optimal solutions.

An alternative approach, proposed in~\cite{Guo-GMACH-22} for feed-forward NNs and in~\cite{Guo-CGMH-23} for LSTMs,  relies on a multi-branch architecture separately processing LoFi and HiFi data, which are then merged by means of a fusion layer combining the information from both fidelities and extracting their correlation to refine the prediction.
In~\cite{Geraci-PGRES-23}, an all-at-once training paradigm is proposed: convolutions are employed in an encoder/decoder architecture to simultaneously compress input data from different fidelities into a latent space, where multi-fidelity data fusion is performed.
This architecture shares conceptual similarities with U-Nets -namely, the encoder/decoder structure with skip connections- which are at the basis of~\cite{Durlofsky-JD-23}.  In this work, a transfer learning procedure is incorporated in a recurrent residual U-Net for which a bulk training is performed using LoFi data, whereas HiFi data are employed to train the output layer and to eventually fine-tune the parameters of the entire network.
A detailed discussion on transfer learning strategies using data generated from a bi-fidelity model is available in~\cite{Doostan-DBRSJD-20}.

Finally,  the paradigm of MFNets~\cite{Gorodetsky-GJGE-20} revists the multi-fidelity setup by defining HiFi and LoFi data, that is, the sources of information of different fidelities, as the nodes of a directed acyclic graph. 
In this context, the connections among nodes represent the correlation existing between the associated data sources and a set of latent variables is used as compact representation of data.
Moreover, in~\cite{Gorodetsky-GJG-21}, data fusion is performed during the all-at-once training that simultaneously uses data from all fidelities leveraging the connections in the previously defined graph.

\begin{remark}[Multi-fidelity surrogates with NN]
In the framework of NN-based multi-fidelity settings,  LoFi data have been generated using ROMs based on both projection approaches and NN architectures.
The paradigm of POD-DL-ROM~\cite{Manzoni-FM-22} first uses randomised POD for dimensionality reduction. Then,  an autoencoder is trained by means of a pre-training using only LoFi data and a fine-tuning relying on HiFi data to further compress the information and construct a latent representation of the coefficients of the reduced basis.
A similar approach exploiting autoencoders for dimensionality reduction and transfer learning for multi-fidelity training is discussed in~\cite{Alonso-SPHXA-24}.
In~\cite{Rozza-DTR-23}, the authors present a multi-fidelity surrogate model blending linear dimensionality reduction (via PODI or gappy POD) and a DeepONet to discover the correlation between the full (HiFi) and the reduced (LoFi) order model and correct the associated mismatch.
An alternative approach to handle multi-fidelity data relies on leveraging the spatial information in the computational meshes of HiFi data by means of graph structures~\cite{Pichi-MPH-24}. In this work,  dimensionality reduction is performed via a NN consisting of feed-forward layers tailored for graph-based data, thus naturally allowing for solutions computed on meshes of different resolution to be employed as input for the training.
Graph-based information is also employed in~\cite{Kutz-KFBK-24} to train a GNN architecture with convolutional layers for dimensionality reduction, following a multi-hierarchical paradigm where HiFi surrogates are trained on the residuals of LoFi ones using transfer learning.
\end{remark}

%---------------------------------------------------------------------
\section{Dataset sampling, enrichment, and augmentation}
\label{sec:Sampling}
%---------------------------------------------------------------------

While multi-fidelity strategies address the efficient exploitation of information available with different levels of accuracy, a complementary approach to achieve efficient surrogate models relies on controlling where and how such information is acquired.
That is, the procedure for the selection of the set of parametric values for which snapshots are computed as full-order instances of the system under analysis plays a critical role in the accuracy and competitiveness of the resulting surrogate model.

The methodologies discussed in this section address complementary aspects of data efficiency in surrogate modelling.
While statistical sampling strategies focus on the selection of parameter instances, adaptive enrichment and active learning aim to iteratively refine the training dataset based on surrogate error indicators, and data augmentation techniques seek to enhance the effective size and diversity of available data.
Although these approaches operate at different stages of the surrogate construction pipeline, they share the common objective of maximising the information content of the training data under limited computational or experimental budgets.
Within this broader perspective, related topics such as sensor placement can be interpreted as experimental design counterparts of parameter-space sampling, and they are particularly relevant when surrogate models are trained from observational data.

Concerning optimal sampling, it is known that Chebyshev points allow to minimise the interpolation error and are thus optimal in dimension 1 (see Remark~\ref{rmrk:Cheb}). 
Nonetheless, in higher dimensional spaces, samples are constructed using a tensorial structure and the number of points increases exponentially with the number of parameters under analysis.
Hence, the goal of sampling is to identify the minimum number of full-order solutions to capture the most relevant features of the problem.
This can be achieved following two main approaches: statistical sampling or adaptive sampling.
On the one hand, statistical sampling employs statistical arguments to select snapshots without exploiting any information on the solution. The goal is to maximise the explored portion of the parametric space. 
These techniques include Latin hypercube sampling~\cite{Conover-79-Latin}, centroidal Voronoi tessellation~\cite{Gunzburger-DFG-99},  and latinised Hammersley and Halton sequences~\cite{Gunzburger-SGB-07}.
On the other hand,  adaptive sampling relies on greedy algorithms to estimate the optimal value of the parameter to be introduced in the dataset, in order to minimise a given measure of the error.
In the following subsections,  some strategies for adaptive sampling, enrichment and augmentation of the dataset of full-order solutions are reviewed.

%---------------------------------------------------------------------
\subsection{Adaptive sampling and enrichment}
\label{sec:Enrichment}
%---------------------------------------------------------------------

Stemming from~\cite{Patera-GP-05,Patera-VP-05}, adaptive sampling strategies have been proposed in ROM literature exploiting the knowledge of the underlying physical model,  e.g.,  via model-constrained optimisation problems~\cite{Willcox-BWG-08} and residual-based error estimates~\cite{Haasdonk-HO-11}. 
In~\cite{Breitkopf-PBBVZ-20},  a strategy to reduce the number of snapshots to be computed was devised in the context of an incremental POD approach sequentially computing snapshots on-the-fly, selecting the values of the parameters identified by an appropriate error estimate.

The above approaches rely on a greedy procedure inspired by the seminal work~\cite{Maday-BMNP-04} and require access to the full-order model employed to compute the snapshots in order to determine the position of the new sampling point,  thus being not suitable for purely data-driven models.
In~\cite{Nadal-NCDFD-15}, a strategy was proposed to select the optimal location of the sampling points using solely a previously computed PGD approximation. The idea consists of locating the new points where the residual of the PGD approximation~\eqref{eq:PGDapprox} is largest, using a criterion based on DEIM.
This consists of selecting, as first sampling point, the point $\bs_1$ in the physical domain where the mode $\fX^1$ reaches a maximum (in absolute value). Thus, the first mode is used to select the first sampling point.
In a similar manner, the second sampling point is defined as the location $\bs_2$ where the absolute value of the function 
\begin{equation}\label{eq:DEIMr2}
\rX^2 :=  \fX^2 - \frac{\fX^2 (\bs_1)}{\fX^1 (\bs_1)}  \fX^1
\end{equation}
achieves a maximum,  with $\rX^2$ being a correction of the second spatial mode $\fX^2$ obtained by subtracting a multiple of the first mode such that $\rX^2(\bs_1) = \bm{0}$.
Similarly, the $k$-th sampling point $\bs_k$ is selected where function $\vert \rX^k \vert$ has a maximum, with
\begin{equation}\label{eq:DEIMrm}
\rX^k :=  \fX^k - \sum_{i=1}^{k-1} d_i^k \fX^i ,
\end{equation}
and the coefficients $d_i^k$ are devised in order to guarantee that $\rX^k(\bs_j)=\bm{0}$ for $j=1,\ldots, k-1$.
Note that the dual vector-function interpretation of modes in the separated description discussed in Remark~\ref{remVectFun} is evoked again here.

This strategy is introduced in the context of sampling the physical space (finding the optimal location of sensors) to properly solve an inverse problem.
Similar ideas are worthy to be employed to optimally sample the parametric space, and producing High-Fidelity snapshots that bring the best information into the surrogate.

Following a similar rationale,  sampling procedures purely data-driven and suitable for non-intrusive coupling with black-box solvers have been devised by many authors. 
In this context,  adaptive sampling strategies need to balance local exploitation and global exploration~\cite{Sobester-SLK-05}. 
The former yields the tendency to place the new sampling point in the region of maximum estimated error, at the risk of exploring only a portion of the parametric space.
The latter is responsible for locating the new sampling point far from previously existing samples (e.g., by means of distance criteria) in order to discover new regions of the parametric space.
CVVor~\cite{Xu-XLWJ-14} combines cross-validation (CV) exploitation to estimate the region of maximum error and Voronoi (Vor) diagrams to perform exploration using a distance-based criterion. 
Similarly,  the smart sampling algorithm (SSA) solves a set of optimisation problems based on CV exploitation and distance-based exploration~\cite{Garud-GKK-17} and the space-filling cross-validation tradeoff (SFCVT)~\cite{Aute-ASAAR-13} relies on leave-one-out CV for exploitation and a distance criterion for exploration.
An alternative approach to drive exploitation is based on an approximated estimate of the gradient of the quantity of interest. In this context,  the local linear approximation (LOLA) couples this strategy with a Voronoi tessellation for exploration~\cite{Crombecq-CGDD-11}, whereas the Taylor expansion-based adaptive design (TEAD) introduced in~\cite{Mo-MLSZYWW-17} performs exploration using a distance-based metric.
In addition, several criteria have been proposed to more systematically balance local exploitation and global exploration. 
For instance, the expected improvement for global fit (EIGF)~\cite{Lam-PhD-08} favours exploitation by selecting the points expected to yield the largest reduction in the global surrogate error, while performing variance-based exploration through uncertainty estimates that guide the sampling toward regions with insufficient information.
The maximising expected prediction error (MEPE)~\cite{Liu-LCO-17} combines exploitation, via the identification of the areas where the surrogate model exhibits high local predictive error, with exploration driven by the maximisation of the expected reduction of the predictive entropy.
The Monte Carlo intersite projection technique (MIPT)~\cite{Crombecq-CLD-11} balances exploitation and exploration by projecting Monte Carlo samples onto locally influential directions to refine regions where the surrogate model shows significant variability, while distributing the samples across intersite locations to capture larger changes in the response surface.

It is worth noting that the techniques described above are mainly tailored to handle sampling when only one model is used to generate data.
In the case of multi-fidelity models, information about fidelity levels is to be accounted for in the adaptive sampling step.
In this context, \cite{Jakeman-JEGG-20} proposes an adaptive version of the MISC algorithm (viz., Section~\ref{sec:MultiInterp}),  to balance the numerical error of the simulations, the error in the surrogate model and the parameter importance sampling, see~\cite{ImportanceSampling-10}, which focuses efforts on critical regions of the parametric space associated with rare events.
In~\cite{Willcox-PKW-18},  the authors employ a hierarchy of surrogate models of different fidelities to perform importance sampling exploiting previously computed biasing densities from LoFi data to accelerate the optimisation process with HiFi data and achieve variance reduction in rare event probability.
An alternative approach based on non-hierarchical fidelities combines accuracy and efficiency in a two-step adaptive sampling strategy by first identifying the HiFi location providing the most potential improvement and selecting the following sample across all fidelities in order to maximise the accuracy improvement normalised by its corresponding computational cost~\cite{Chen-CWCC-24}.
Finally,  the mixed adaptive sampling algorithm (MASA) proposed in~\cite{Cremaschi-EC-14} uses a query-by-committee approach for local exploitation in order to select the new candidate sample for which the discrepancy among the predictions performed by the different surrogates (or fidelities) under analysis is the most significant,  while maximising the sample distance for exploration purposes.
For a detailed discussion on adaptive sampling strategies for surrogate models,  interested readers are referred to~\cite{Liu-LPC-18,Fuhg-FFN-21,Mainini-dFNM-24}.

%---------------------------------------------------------------------
\subsection{Data augmentation}
\label{sec:Augmentation}
%---------------------------------------------------------------------

As previously mentioned, the quality of the constructed surrogate models significantly depends upon the richness of the dataset employed for training.
Nonetheless, the construction of an appropriate training set for the surrogate model can become computationally unfeasible when the full-order models are nonlinear,  feature a large number of parameters, and the cost of each snapshots is particularly expensive (e.g., for large-scale, coupled, multi-physics, or multi-scale problems).

To circumvent this issue,  data augmentation strategies aim to reduce the computational investment for the construction of the initial dataset by engineering artificial snapshots providing meaningful additional information, not present in the training set.
The idea is to generate new samples without evaluating the full-order model in new parametric values.
Different approaches have been proposed in the literature in recent years.
In~\cite{diez2021nonlinear},  a set of inexpensive, algebraic operations among snapshots were performed to enrich the information of the training set to construct a POD approximation.
In~\cite{Muixi-MZGD-25}, the authors propose to design new, artificial snapshots compliant with the underlying physics by enforcing conservation of mass and conservation of momentum in the data augmentation procedure.
In~\cite{Nguyen-25}, a generative approach exploiting nonlinear activation functions is used to combine the information of previously existing snapshots.

The resulting data augmentation procedures are able to emancipate the quality of the a posteriori ROM approximation from the quality of the initial dataset. This is achieved by incorporating new information in uniformly-sampled, and possibly poorly populated, existing datasets.
It is worth noting that this enrichment strategy is safe to use if the resulting dataset is used to solve the equation as in POD-RB, because physics is thus enforced in the final solution when solving the equations in the reduced space. 
For other ROM strategies, for instance PODI, the final result relies on the physical pertinence of the augmented data. Thus, in these cases (very important if the methodology is non-intrusive), the method used to produce the new, artificial snapshots has to be credible, and appropriate techniques to assess their reliability are required.

%---------------------------------------------------------------------
\section{Conclusions and future perspectives}
\label{sec:Conclusion}
%---------------------------------------------------------------------

This article presented an overview of established techniques to construct (parametric) functional surrogates of complex engineering systems using physics-based and data-driven approaches.
These strategies are becoming increasingly relevant across a wide range of application domains requiring efficient computational workflows for repeated evaluations of complex models.
In this context, this work has revisited surrogate modelling from a functional approximation perspective, highlighting common methodological foundations across different communities and discussing representative techniques for dimensionality reduction and functional representation.

While this review focuses on the methodological aspects of surrogate modelling of parametric systems, the field is rapidly evolving, beyond academic research.
In particular, surrogate models are attracting growing attention from industry due to their potential to enable efficient analysis of complex systems in data-rich environments such as smart cities, manufacturing 4.0, personalised healthcare, and sustainability, where real-time analysis is often required.

%---------------------------------------------------------------------
\subsection{Critical discussion}
%---------------------------------------------------------------------

Following the functional approximation viewpoint adopted in this review, building a surrogate model can be interpreted as the combination of an approximation strategy with the construction of a suitable reduced representation of the solution manifold. The methods presented in Sections~\ref{sec:POD}, \ref{sec:PGD}, \ref{sec:NN}, and~\ref{sec:MultiFidelity} primarily differ in how these two aspects are treated but several common considerations emerge. The following discussion synthesises these insights through the criteria introduced in Section~\ref{sec:Intro}, emphasising the main methodological trade-offs in terms of intrusiveness, data requirements, computational efficiency, expressiveness, and interpretability.

Projection-based approaches rely on the construction of low-dimensional linear subspaces that approximate the dominant features of the solution manifold. When the governing equations are enforced at the reduced level (e.g., in POD-RB), the resulting surrogate models are interpretable but often require intrusive access to the full-order solvers limiting their applicability in industrial settings. Although their mathematical foundations allow for rigorous error estimation and stability analysis, the reliance on linear subspaces may limit their expressiveness for strongly nonlinear problems or for systems exhibiting complex parametric dependence.

Separated representations such as PGD enable explicit parametric dependence of the surrogate representation, circumventing the need for the solution of any reduced problem during the online evaluation phase. This can significantly alleviate the curse of dimensionality in problems with many parameters. PGD thus offers a powerful framework for parametric modelling when the solution admits a low-rank separable structure. Nevertheless, the construction of such separated representations can become challenging for highly nonlinear or weakly separable systems, and the efficiency of the method may depend strongly on the structure of the underlying problem.

In contrast, purely data-driven approaches (e.g., PODI, a posteriori PGD, NNs, $\ldots$) offer flexible functional approximations capable of capturing complex nonlinear relations directly from data, but they typically provide less explicit physical interpretability. Their ability to approximate highly nonlinear functions makes them attractive for problems where classical linear reduced representations are insufficient. Among these strategies, NNs are particularly appealing because they can be used not only to approximate parametric mappings but also to learn nonlinear reduced representations of the solution manifold itself. Performance of these models typically depends on the availability of sufficiently rich training datasets, making them well-suited for data-rich environments and less appropriate for data-scarce settings typically arising in computational science and engineering.

Finally, multi-fidelity strategies aim to combine information from models with different levels of accuracy in order to reduce the computational cost of surrogate construction. By exploiting correlations between low- and high-fidelity models, these approaches enable efficient use of limited high-fidelity data and provide a flexible framework to integrate heterogeneous sources of information. Their effectiveness, however, depends on the existence of suitable hierarchies between the fidelities, and on the strength of their correlations. In this context, it is thus crucial to establish a trade-off between discretisation and model errors.

To summarise, projection-based and physics-informed approaches provide structure, interpretability, and strong theoretical foundations, whereas data-driven approaches offer greater flexibility to capture complex nonlinear behaviours. Multi-fidelity frameworks complement these methodologies by improving data efficiency when multiple information sources are available. Rather than competing alternatives, these strategies should therefore be viewed as complementary tools that can be combined to address different modelling challenges in parametric systems.

These differences highlight how the choice of surrogate modelling strategy ultimately depends on the characteristics of the underlying physical model, the discretisation, the dimensionality of the parameter space, and the availability of training data. As a consequence, direct quantitative comparisons between the different approaches are inherently problem-dependent. 
Hence, a systematic benchmark assessment would require carefully designed test cases and consistent evaluation settings tailored to specific applications and therefore lie beyond the scope of the present methodological review.

%---------------------------------------------------------------------
\subsection{Open challenges}
%---------------------------------------------------------------------

Building upon the methodological discussion presented throughout the manuscript, several aspects remain to be further investigated to improve both the understanding and the practical deployment of surrogate modelling techniques. The following subsections present a non-exhaustive list of open problems and future perspectives in the field, spanning different levels of maturity and research readiness. 
These challenges are organised into three categories:
\begin{itemize}
\item scalability and computational efficiency of the surrogate modelling techniques;
\item accuracy, reliability, and robustness of surrogate predictions;
\item expressive capabilities of surrogate models to describe the underlying relations between user-defined input parameters and output quantities of interest.
\end{itemize}
Some of the topics have already been addressed in the literature, although the proposed solutions are not yet fully mature or production-ready.  Others represent active research directions currently pursued by several groups worldwide. Finally, we also include a set of emerging topics that may have a significant impact on the development of surrogate modelling methodologies in the coming years.

%---------------------------------------------------------------------
\subsubsection{Scalability and computational efficiency}
%---------------------------------------------------------------------

A central objective for the future development of surrogate modelling techniques is their deployment in realistic, large-scale engineering applications with many parameters. Achieving this vision requires methodologies that remain efficient and scalable when confronted with high-dimensional parameter spaces, complex physical models, and that can fully exploit modern and portable computing infrastructures.
\begin{itemize}
\item \textbf{High parametric dimensionality -} Many realistic engineering and scientific applications feature an extreme number of parameters.  Dimensionality reduction techniques suitable to handle such high-dimensional problems will thus play a crucial role in developing, deploying, and adopting surrogate models for realistic applications. Besides the methodologies reviewed in this work, alternative successful dimensionality reduction approaches have been proposed in recent years, including low rank approximations~\cite{Markovsky-book-12}, tensor-train decompositions~\cite{Oseledets-11}, and dynamic mode decomposition~\cite{Kutz-TRLBK-14}. 
\item \textbf{Large physical scale -} The construction of physics-based and data-driven surrogates rely on the solution of multiple instances of a, possibly large-scale, full-order problem. To make such computations affordable, several aspects require further development, including the construction of local ROMs via domain decomposition~\cite{Ohlberger-BIORSS-20,Barnett:2022:Sandia,MG-DEG-24,MG-DEG-26} or dynamic substructuring~\cite{Park-PP-20,Chatzi-VTABQC-22,Berthet-BPCG-24,Sengupta-SC-25} and the interface with high-performance computing infrastructures, see, e.g., ~\cite{Lehmkhul-EMBVRL-25}. 
These solutions offer promising directions to exploit locality, modularity, and parallelism in large, multi-component systems, guaranteeing scalability and computational efficiency as model size increases.
\item \textbf{Computing paradigms -} To achieve democratisation and upscaling of digital twin technologies, the reviewed strategies to construct surrogate models need to be extended in order to account for available computing paradigms. On the one hand, many industries rely on legacy systems (e.g., Fortran-based as discussed in~\cite{FortranLegacy-22}) and computing infrastructures lacking the connectivity to be seamlessly integrated with interactive digital twins according to the Industry 4.0 paradigm. On the other hand, the upfront investment in hardware and software for access to (and maintenance of) parametric surrogate modelling solutions might represent a bottleneck for small and medium-sized enterprises. To tackle the former, tailored wrappers of surrogate methodologies need to be developed in upcoming years. For the latter,  models based on cloud computing and digital twin as a service (DTaaS)~\cite{DTaaS-21} are expected to reduce entry barriers for smaller organisations. Moreover, applications requiring real-time or nearly real-time feedback might encounter a practical bottleneck in the communication latency between sensors and the surrogate model itself, limiting its capability to provide instantaneous responses. In this context,  edge computing offers under-explored possibilities for digital twins by performing computations locally on devices, reducing data interchange with cloud servers and faster updates~\cite{EdgeDT-24}. Finally,  new computing paradigms, such as quantum computing, are to be further explored to assess their suitability to construct surrogate models of complex systems, as discussed in a recent work which employs quantum unconstrained binary optimisation and quantum approximation optimisation algorithm for dynamic mode decomposition~\cite{QuantumROM-24}.
\end{itemize}

%---------------------------------------------------------------------
\subsubsection{Accuracy, reliability, and robustness}
%---------------------------------------------------------------------

As surrogate models become increasingly integrated into decision-making workflows, ensuring the reliability and robustness of their predictions becomes a fundamental requirement. This calls for continued advances in data management, validation and standardisation practices, and modelling strategies capable of supporting trustworthy decision-making under uncertainty.
\begin{itemize}
\item \textbf{Data quality -} The accuracy and robustness of surrogate models significantly depend upon the quality and fidelity of input data used for training. The multi-fidelity framework discussed in Section~\ref{sec:MultiFidelity}, coupled with tailored uncertainty quantification strategies, is expected to further develop in order to account for challenging problems such as noisy data~\cite{Diez-PWBSVD-23,Kent-KTGH-26}, non-hierarchical fidelities~\cite{Willcox-LAW-15,Zhang-ZWJCZ-22},  missing and incomplete data. Moreover, an alternative non-recursive framework for multi-fidelity has been recently proposed to circumvent accuracy limitations of existing nested sampling strategies by leveraging all correlations among the different sources of information~\cite{Gorodetsky-GGEJ-20}.
\item \textbf{Verification \& Validation -} Reproducibility of results is pivotal for scientific advancement. Nonetheless, in the context of scientific machine learning,  lack of standardised benchmarks and limited access to data and source codes have limited the potential impact of performed research. To address this issue, a recent manuscript sketches a roadmap towards trustworthy scientific machine learning, leveraging existing knowledge from verification \& validation in the computational science and engineering community~\cite{Jakeman-JBMO-Preprint-25}. Additional key aspects to be accounted for in upcoming years include findability, accessibility, interoperability, and reusability of research data generated in the field according to the FAIR principles~\cite{FAIR-16}.
\item \textbf{Calibration \& uncertainty quantification -} A critical challenge for the reliable deployment of surrogate models in digital twin applications is the calibration of predictive uncertainty, especially since classical built-in uncertainty expressions are not necessarily well calibrated in practice, see, e.g.,~\cite{Bastos-BOH-09} for the case of GPR. This issue is of particular relevance in data-scarce and extrapolative regimes~\cite{Forrester-FSK-08,Maday-BFM-23}, as the availability of reliable and well-calibrated uncertainty estimates is a necessary condition to enable trustworthy solutions for decision-making under uncertainty~\cite{Lin-LBD-21,Thelen-TZFLGYTMHH-23,Willcox-SHCOBPWCJ-25}.
\end{itemize}

%---------------------------------------------------------------------
\subsubsection{Model capabilities and expressiveness}
%---------------------------------------------------------------------

Beyond efficiency and reliability, future surrogate modelling frameworks will need to provide richer modelling capabilities in order to capture increasingly complex physical phenomena and adapt to evolving and heterogeneous data streams. Advancing the expressive power of surrogate models while preserving interpretability and generalisation represents a key research direction for the coming years.
\begin{itemize}
\item \textbf{Model complexity -} Nonlinear problems, multi-scale phenomena, chaotic dynamics, discontinuous and bifurcating solutions still entail significant computational challenges for existing surrogate modelling methods and will require dedicated efforts in the future. In this context, alternative approaches based on model discovery, such as sparse identification of nonlinear dynamics (SINDy)~\cite{Brunton-BPK-16}, have been proposed to circumvent the computational difficulties of approximating complex physical models during surrogate construction.
\item \textbf{Generalisation -} Availability of \emph{sufficiently rich} datasets for training is crucial to construct an accurate surrogate model with generalisation capabilities, especially concerning parameter extrapolation and domain-awareness. Attempts to tackle this problem are currently being explored in the literature using a variety of solutions including importance sampling to identify critical regions of the parametric space~\cite{Peherstorfer-PCMW-16}, generation of artificial snapshots by means of data augmentation~\cite{Muixi-MZGD-25} and optimal transport~\cite{Torregrosa-TCAHC-22}, transfer learning~\cite{Kadeethum-KOCVY-24} and meta-learning~\cite{Beer-CDFBY-22} to leverage prior knowledge without significant retraining, and generative models such as variational autoencoders~\cite{Vinuesa-WSVV-24} and generative adversarial networks~\cite{Rozza-CDR-24}. 
\item \textbf{Self adaptation -} Developing models that can dynamically adapt themselves to new problems and new data is crucial to minimise retraining costs.  In this context,  strategies from \emph{online learning} and \emph{continual learning} are expected to provide significant advances in upcoming years. On the one hand, online learning focuses on adapting, in real-time, a model designed for one task by accounting for new streaming data. Proposed solutions in this context include dynamical low-rank approximation~\cite{Lubich-KL-07},  low-rank updates in dynamic data-driven application systems (DDDAS)~\cite{Peherstorfer-PW-15},  on-the-fly ROMs~\cite{Breitkopf-PBBVZ-20,Babaee-RNB-21}, online adaptive basis refinement~\cite{Carlberg-EC-20}, dynamic mode decomposition~\cite{Kutz-TRLBK-14},  sequential data assimilation and integration using Kalman filtering~\cite{Poux-MP-17} and transport maps~\cite{Chamoin-Preprint-24}. On the other hand, continual learning~\cite{Kirkpatrick-KPRVDRMQRGHCKH-17} aims to dynamically learn new tasks while retaining prior knowledge acquired during training, thus offering potentially enormous advantages for self-adapting surrogate models enabling the adoption of digital twins in complex, time-dependent processes.  Similarly, the paradigm of reinforcement twinning~\cite{Mendez-SMPAvdBM-24} has been recently proposed to exploit reinforcement learning for digital twins construction.
\item \textbf{Interpretability \& Explainability -} Whilst non-intrusiveness represents a major advantage for the applicability of reduced order methodologies to realistic scenarios,  relying exclusively on \emph{black-box} surrogates makes it difficult to understand and justify the resulting predictions. In this context,  increasing attention is being attracted by interpretability~\cite{Ramabathiran-RR-21} of the relation between inputs and outputs of the surrogates in terms of existing scientific knowledge, and explainability~\cite{Vinuesa-CHDQLLMHLMV-24} of the models in terms of human-understandable concepts.
\end{itemize}
%---------------------------------------------------------------------

%---------------------------------------------------------------------
\section*{Acknowledgements}
%---------------------------------------------------------------------
This work was partially supported by the Spanish Ministry of Science, Innovation and Universities and the Spanish State Research Agency MICIU/AEI/10.13039/501100011033 (Grant agreement No. TED2021-132021B-I00 to MG; PID2023-149979OB-I00 to MG; PID2023-153082OB-I00 to PD), EIT KIC RawMaterials Project Agreement 21123,  European Commission MSCA Ref. 101120556-EARTHSAFE, and the Generalitat de Catalunya (Grant agreement No. 2021-SGR-01049). MG is Fellow of the Serra H\'unter Programme of the Generalitat de Catalunya.

%---------------------------------------------------------------------
\section*{Conflict of interest}
%---------------------------------------------------------------------
The authors declare that they have no conflict of interest.

%---------------------------------------------------------------------
\bibliographystyle{plain}
\bibliography{Ref}
%---------------------------------------------------------------------

\end{document}